\documentclass[a4paper,11pt]{amsart}
\usepackage{amsmath}
\usepackage[T1]{fontenc}
\usepackage{amssymb}
\usepackage{amsthm}

\newtheorem{thm}{Theorem}[section]

\newtheorem{theorem}[thm]{Theorem}

\newtheorem{lemma}[thm]{Lemma}
\newtheorem{proposition}[thm]{Proposition}

\newtheorem{rem}{Remark}[section]

\numberwithin{equation}{section}

\title [ Sharp Well-posedness  for  (KPBII)  equation]
{ Sharp  well-posedness  for Kadomtsev-Petviashvili-Burgers (KPBII) equation in $\mathbb R^2$}

\author{Bassam Kojok}
\date{\today}
\subjclass[]{}
\address{Laboratoire Analyse, G\'eom\'etrie et Applications, Institut Galil\'ee, Universit\'e Paris-Nord, 93430 Villetaneuse, France.}
\email{ kojok@math.univ-paris13.fr}
\begin{document}
\maketitle
\begin{abstract}
We prove global  well-posedness  for the  Cauchy problem associated with
the Kadomotsev-Petviashvili-Burgers equation (KPBII) in $\mathbb
R^2$ when the initial value  belongs to the anisotropic Sobolev space
$H^{s_1,s_2}(\mathbb R^2)$ for all $s_1>-\frac12$ and $s_2\geq0$. On
the other hand, we prove in some sense that our result is sharp.
\end{abstract}

\section{Introduction}
We shall study the initial value problem of the
Kadomtsev-Petviashvili-Burgers (KPBII) in $\mathbb R^2$~:
\begin{equation}\label{KPB}
\left\{\begin{array}{l}
\left(\partial_tu+u_{xxx}-u_{xx}+uu_x\right)_x+u_{yy}=0,\\
u(0,x,y)=\varphi(x,y).
\end{array}\right.\end{equation}
This equation is a dissipative version of the
Kadomtsev-Petviashvili-II equation (KPII)~:
\begin{equation}\label{KP}
\left\{\begin{array}{l}
\left(\partial_tu+u_{xxx}+uu_x\right)_x+u_{yy}=0,\\
u(0,x,y)=\varphi(x,y).
\end{array}\right.\end{equation}

 The (KP) equation is a universal model for nearly
one directional weakly nonlinear dispersive waves with weak
transverse effects. It  is  a natural two
dimensional extension of the celebrated (KdV) equation~:
$$u_t+u_{xxx}+uu_x=0.$$
In some typical situations, it is not possible to neglect dissipative
effects (due to viscosity effects in magneto sonic waves damped by
electron-ion collisions for example), and this can  lead to the
KdV-Burgers equation (cf. \cite{OtSu})~:
$$
\partial_tu+u_{xxx}+uu_x-u_{xx}=0.
$$
It  is then widely accepted that the (KPB) equation is a  natural
model for the propagation of the two dimensional damped waves. Note
that as we are interested in  nearly one directional propagation,
the dissipative term only acts in the main direction of propagation
in (\ref{KPB}).\\
Bourgain had developped a new method, clarified by Ginibre in
\cite{Ge}, for the study of Cauchy problem associated with
dispersive non-linear equations. This method was successfully
applied to  Schr\"odinger, (KdV) as well as (KPII) equations. It
was  shown by Molinet-Ribaud \cite{MoRi}  that the Bourgain spaces can
be used to study the  Cauchy problems associated to semi-linear equations
with a linear part containing both dispersive and dissipative
terms (and consequently this applies to (KPB) equations).\\
For the Cauchy problem associated to (KPII) equation, the local
existence is proved by Bourgain \cite{Bo93}    when the initial value
is in the space $L^2(\mathbb R^2)$ and by Takaoka-Tzvetkov
\cite{TaTz}  when the initial value $\varphi\in H^{s_1,s_2}(\mathbb
R^2)$ with $s_1>-\frac13$ and $s_2\geq0$.\\
By introducing a Bourgain space associated to the usual (KPII)
equation (related only to the dispersive part of the linear symbol
of (\ref{KPB})), Molinet-Ribaud \cite{MoRi} had proved global
existence for  the Cauchy problem associated to the (KPBII) equation
when the initial value is in $L^2(\mathbb R^2)$.\\
In this paper, we prove local existence for (\ref{KPB}) with initial
value $\varphi\in H^{s_1,s_2}(\mathbb R^2)$ when  $s_1>-\frac12$ and
$s_2\geq0$. Following \cite{MoRi2}, we introduce a Bourgain space associated to the (KPBII)
equation. This space is in fact the intersection of the space
introduced in \cite{Bo93} and of a Sobolev  space. The advantage of
this space is that it contains both the dissipative and dispersive
parts of the linear symbol  of (\ref{KPB}). \\
We prove also that our local existence theorem is optimal by
constructing a counter example showing that the application
$\varphi\mapsto u$ from  $H^{s_1,s_2}$ to $C([0,T];H^{s_1,s_2})$ can not be
 regular for $s_1<-\frac12$ and $s_2=0$.\\
The paper is organized as follows. In Section \ref{Not}, we
introduce our notations and we give an extension of the semi-group
of the (KPBII) equation by a linear operator defined on all the real
axis. In Section \ref{EstLin} we derive  linear estimates and some smoothing
 properties for the operator $L$ defined by (\ref{ForcingOp2}) in the
Bourgain spaces . In Section \ref{EstStr} we
state Strichartz type estimates for the (KP) equation which yield
bilinear estimates in Section \ref{EstBL}. In Section \ref{RegEL},
 using bilinear estimates,  a standard fixed point
argument and some smoothing properties, we prove uniqueness and global existence of the solution of
(\ref{KPB}) in anisotropic sobolev space $H^{s_1,s_2}(\mathbb R^2)$
with $s_1>-\frac12$ and $s_2\geq0$. Finally,  we construct  in Section \ref{ConExp}   a sequence of initial values which ensures
that our local existence result is optimal if one requires the  smoothness of the flow-map. Note that there is no scaling for (\ref{KPB}) and that, on the other hand, $H^{-1/2,0}$ is critical for the scaling of (\ref{KP}).
\section{Notations and main results}\label{Not}
We will use $C$ to denote various time independent constants, usually depending only upon $s$.
In case a constant depends upon other quantities, we will try to make it explicit.
 We use $A \lesssim B$ to denote an estimate of the form $A\leq C B$. similarly,
 we will write $A\sim B$ to mean $A \lesssim B$  and $B \lesssim A$. We  writre
 $\langle \cdot\rangle :=(1+|\cdot|^2)^{1/2}\sim 1+|\cdot| $. The notation $a^+$ denotes
$a+\epsilon$ for an arbitrarily small $\epsilon$. Similarly $a-$ denotes $a-\epsilon$.
 For $b\in\mathbb R$, we denote respectively by $H^b(\mathbb R)$  and ${\dot{H}}^b(\mathbb R)$
the  nonhomogeneous and homogeneous Sobolev spaces  which are endowed with the  following norms~:
\begin{equation}\label{N1}
||u||^2 _{H^b}=\int_{\mathbb R}{\langle\tau\rangle}^{2b}|\hat u(\tau)|^2 d\tau,\quad
||u||^2 _{{\dot H}^b}=\int_{\mathbb R} |\tau|^{2b} |\hat u(\tau)|^2 d\tau
\end{equation}
where $\hat .$ denotes the Fourier transform from $\mathcal S'(\mathbb
R^2)$ to $\mathcal S'(\mathbb R^2)$ which is defined by~:
$$\hat f(\xi):=\mathcal F(f)(\xi)=\int_{\mathbb R^2}e^{i\langle
\lambda,\xi\rangle}f(\lambda)d\lambda,\quad \forall f\in \mathcal
S'(\mathbb R^2).$$
Moreover, we  introduce the corresponding space (resp space-time)  Sobolev spaces   $H^{s_1,s_2}$
(resp  $ H^{b,s_1,s_2}$)  which are defined by~:
\begin{equation}\label{Esp1}
H^{s_1,s_2}(\mathbb R^2)=:\{ u \in \mathcal S^{'}(\mathbb R^2); ||u||_{H^{s_1,s_2}}(\mathbb R^2)
<+\infty \},
\end{equation}
\begin{equation}\label{Esp2}
H^{b,s_1,s_2}(\mathbb R^2)=:\{ u \in \mathcal S^{'}(\mathbb R^3); ||u||_{H^{b,s_1,s_2}}(\mathbb R^3)
<+\infty \}
\end{equation}
where,
\begin{equation}\label{Ns}
||u||^2 _{H^{s_1,s_2}}=\int_{\mathbb R^2}{\langle\xi\rangle}^{2s_1}{\langle\eta\rangle}^{2s_2}
|\hat u(\nu)|^2 d\nu,
\end{equation}
\begin{equation}\label{Nbs}
||u||^2 _{H^{b,s_1,s_2}}=\int_{\mathbb R^2}{\langle\tau\rangle}^{b}{\langle\xi\rangle}^{2s_1}
{\langle\eta\rangle}^{2s_2}|\hat u(\tau,\nu)|^2 d\nu d\tau,
\end{equation}
and $\nu =(\xi,\eta)$. Let $U(\cdot)$ be the unitary group in $H^{s_1,s_2}$, $s_1$, $s_2\in\mathbb R$, defining the free evolution of the (KP-II) equation, which is given by
\begin{equation}\label{Gr}
U(t)=\exp(itP(D_x,D_y)),
\end{equation}
 where $P(D_x,D_y)$ is the Fourier multiplier with symbol $P(\xi,\eta)=\xi^3 - \eta^2/ \xi$.
 By the  Fourier transform, (\ref{Gr}) can be written like~:
\begin{equation}\label{Grf}
\mathcal{F}_x(U(t)\phi)=\exp(itP(\xi,\eta))\hat{\phi}, \quad
\forall\phi \in \mathcal S^{'}(\mathbb R^2),\quad t\in \mathbb R.
\end{equation}
Also, by the Fourier transform, the linear part of the  equation \ref{KPB} can be written as~:
\begin{equation}\label{Symb}
i(\tau -\xi^3 -\eta^2/\xi) + \xi^2=:i(\tau - P(\eta,\xi)) + \xi^2.
\end{equation}
This leads us, as in  \cite{MoRi},  to introduce a  Bourgain space  which is in relation
with  both the  dissipative and  dispersive  parts  of (\ref{KPB}) at the same time, we define this space by~
\begin{equation}\label{EspBo}
X^{b,s_1,s_2}=\{ u\in \mathcal S^{'}(\mathbb R^3), ||u||_{X^{b,s_1,s_2}}<\infty \}
\end{equation}
 equipped with the norm
\begin{equation}\label{NoBo}
||u||_{X^{b,s_1,s_2}}=||\langle i\sigma + \xi^2 \rangle^b \langle \xi\rangle^{s_1}
\langle \eta \rangle ^{s_2} \hat w(\tau,\nu)||_{L^2(\mathbb R^3)}
\end{equation}
where,\\
 $\sigma=\tau -P(\nu),\quad \nu=(\xi,\eta)\in \mathbb R^2.$

\begin{rem}\label{RqEs}
It is worth noticing that $X^{b,s_1,s_2}$ is the intersection of the Bourgain  space associated with  the dispersive part of equation (\ref{KPB})  and Sobolev  space. Indeed, by noticing that
 $\mathcal F(U(-t)u)(\tau,\nu)=\mathcal F(u)(\tau + P(\nu),\nu)$ and by the change of variable
$\tau \longrightarrow \tau - P(\nu)$, one sees that
 \begin{eqnarray*}
||u||_{X^{b,s_1,s_2}}&=&||\langle i\tau + \xi^2 \rangle^b \langle \xi\rangle^{s_1}
\langle \eta \rangle ^{s_2} \hat w(\tau +P(\nu),\nu)||_{L_{\tau,\nu}^2(\mathbb R^3)}\\
&=&||\langle i\tau + \xi^2 \rangle^b \langle \xi\rangle^{s_1}
\langle \eta \rangle ^{s_2}\mathcal F(U(-t)u)(\tau,\nu) ||_{L_{\tau,\nu}^2(\mathbb R^3)}\\
&\sim &||\langle \tau   \rangle^b \langle \xi\rangle^{s_1}
\langle \eta \rangle ^{s_2}\mathcal F(U(-t)u)(\tau,\nu) ||_{L_{\tau,\nu}^2(\mathbb R^3)}\\
&+&| \langle \xi\rangle^{s_1+2b}
\langle \eta \rangle ^{s_2}\mathcal F(U(-t)u)(\tau,\nu) ||_{L_{\tau,\nu}^2(\mathbb R^3)}\\
&=&||U(-t)u||_{H^{b,s_1,s_2}} + ||u||_{L^2_t H^{s_1+2b,s_2}}.
\end{eqnarray*}\hfill$\Box$
\end{rem}
For $T>0$, we define the restricted spaces  $X^{b,s_1,s_2}_T$ by the norm
\begin{equation}\label{EspBoLocal}
||u||_{X_T^{b,s_1,s_2}}=\inf_{w\in X^{b,s_1,s_2}}\{||w||_{X^{b,s_1,s_2}}; w(t)=u(t) \text{ on }
  [0,T]\}.
\end{equation}
We denote  by $W(\cdot)$ the semigroup associated with the free evolution of (\ref{KPB}),
\begin{equation}\label{SemiGr}
\mathcal{F}_x(W(t)\phi)=\exp(itP(\xi,\eta)-|\xi|^2t)\hat{\phi}, \quad
\forall\phi \in \mathcal S^{'}(\mathbb R^2),\quad t\geq 0.
\end{equation}
Also, we can extend $W$ to a linear operator defined on the whole real axis by setting,
\begin{equation}\label{SemiGrPr}
\mathcal{F}_x(W(t)\phi)=\exp(itP(\xi,\eta)-|\xi|^2|t|)\hat{\phi}, \quad
\forall\phi \in \mathcal S^{'}(\mathbb R^2),\quad t\in \mathbb R .
\end{equation}
By the  Duhamel integral formulation, the equation \ref{KPB} can be written,
\begin{equation}\label{FormDuh}
u(t)=W(t)\phi - \frac12 \int_0^t W(t-t')\partial_x(u^2(t'))dt', \quad t\geq 0.
\end{equation}
To prove the local existence result, we will apply a fixed  point argument
  to a truncated version of (\ref{FormDuh})  , which is defined on all the real axis by
\begin{equation}\label{FormDuhTr}
u(t)=\psi(t)[ W(t)\phi - \frac{\chi_{\mathbb R_+}(t)}{2} \int_0^t W(t-t')
\partial_x(\psi^2_T(t')u^2(t'))dt'],
\end{equation}
where $t\in \mathbb R$ and $\psi$ indicates a time cutoff fonction~:
\begin{equation}\label{Cutoff}
\psi \in C_0^{\infty}(\mathbb R),\quad \text{sup }\psi\subset [-2,2], \quad \psi=1\text{ on }[-1,1],
\end{equation}
and $ \psi_T(.)=\psi(./T).$
\begin{rem}\label{RemDequiv}
It is clear that if $u$ solves (\ref{FormDuhTr}) then $u$ is a solution of (\ref{FormDuh}) on $[0,T]$,
 $T<1$. Thus it is sufficient to solve (\ref{FormDuhTr})  for a small time ($T<1$ is enough).
\end{rem}
Let us now state our results:
 \begin{theorem}\label{Th1}
Let $s_1>-1/2$, $s_2\geq 0$, $s_c^1\in]-1/2,\min(0,s_1)]$  and $\phi \in H^{s_1,s_2}$. Then for
 $s_1\geq s_c^1$  there exists a time
 $T=T(||\phi||_{H^{s_c^1,0}})>0$ and a unique solution $u$ of  (\ref{KPB}) in
\begin{equation}\label{Th1Esp}
Y_T=C([0,T];H^{s_1,s_2})\cap X_T^{1/2,s_1,s_2}
\end{equation}
Moreover, $u\in C(\mathbb R_+;H^{s_1,s_2})\cap C(\mathbb R_+^*;H^{\infty,s_2})$ and  the map $\phi\longmapsto u$ is $C^\infty$ from $H^{s_1,s_2}$ to $Y_T$.
$\hfill{\Box}$
\end{theorem}

\begin{theorem}\label{Th2}
Let $s<-1/2$. Then it  does not exist a time $T>0$  such that the equation (\ref{KPB}) admits a unique
solution in $C([0,T[,H^{s,0})$  for any  initial data in some ball of  $H^{s,0}(\mathbb R^2)$
 centered at the origin and such that the map
\begin{equation}\label{Th2Flow}
\phi\longmapsto u
\end{equation}
is $C^2$-differentiable  at the origin from $H^{s,0}$ to $C([0,T],H^{s,0})$.
$\hfill{\Box}$
\end{theorem}
\section{Linear estimates in $X^{b,s_1,s_2}$}\label{EstLin}
 In this section we study  both the  free and the forcing  terms of the integral equation (\ref{FormDuhTr})
to obtain certain estimates necessary to apply a  fixed point argument. The results of this section
 are essentially contained in \cite{MoRi}.The following lemma will
be of constant use in this section~:
\begin{lemma}\label{LemConst}
 let $b\in \mathbb R$ and $\lambda>0$. Then
\begin{equation}\label{LemH}
||f(\lambda t)||_{H^b}=(\lambda^{-1/2} + \lambda^{b-1/2})||f(t)||_{H^b},
\end{equation}
\begin{equation}\label{LemHPo}
||f(\lambda t)||_{\dot H^b}= \lambda^{b-1/2}||f(t)||_{\dot H^b}.
\end{equation}$\hfill{\Box}$
\end{lemma}
\begin{proposition}\label{ProFree}
Let $s_1$, $s_2\in\mathbb R$ and $0\leq b\leq 1/2$. For all $f\in H^{s_1,s_2}$ we have,
\begin{equation}\label{FreeEst}
||\psi(t)W(t)\phi||_{X^{b,s_1,s_2}}\leq C ||\phi||_{H^{s_1 +2b-1,s_2}}
\end{equation}
$\hfill{\Box}$
\end{proposition}

 $\bold{Proof}$. By definition of $W(\cdot)$ and $X^{b,s_1,s_2}$, and Using the change of variable
 $\tau\longmapsto\sigma:=\tau -P(\nu)$ we have,
\begin{eqnarray}\label{ProFreeDec}
\nonumber||\psi(t)W(t)\phi||_{X^{b,s_1,s_2}}&=&||\langle i\sigma + \xi^2 \rangle^b \langle
\xi\rangle^{s_1}\langle \eta \rangle ^{s_2} \mathcal{F}_t(\psi(t)e^{-|t|\xi^2}
e^{itP(\nu)})(\tau)||_{L^2}\\
\nonumber&=&||\langle i\tau + \xi^2 \rangle^b \langle
\xi\rangle^{s_1}\langle \eta \rangle ^{s_2} \mathcal{F}_t(\psi(t)e^{-|t|\xi^2})(\tau)||_{L^2}\\
\nonumber&\leq&|| \langle \xi\rangle^{s_1+2b}\langle \eta \rangle ^{s_2}
|| \mathcal{F}_t(\psi(t)e^{-|t|\xi^2})(\tau)||_{L^2_\tau}||_{L^2_\nu}\\
&+&|| \langle \xi\rangle^{s_1}\langle \eta \rangle ^{s_2}||\langle\tau\rangle^b
\mathcal{F}_t(\psi(t)e^{-|t|\xi^2})(\tau)||_{L^2_\tau}||_{L^2_\nu}
\end{eqnarray}
For $\xi$ fixed, let $0\leq b\leq 1/2$,  we  pose $G_\xi(\tau)=\langle\tau\rangle^b
\mathcal{F}_t(\psi(t)e^{-|t|\xi^2})(\tau)$. we notice, that as in  \cite{MoRi}, we have the
following estimate~:
\begin{equation}\label{FreeMoRi}
||G_\xi||_{L^2_\tau}(\mathbb R)\leq C\langle\xi\rangle^{2b-1},\quad \forall 0\leq b\leq1/2
\end{equation}
By combining these two last inequalities, one obtains the  requiered  result.
$\hfill{\Box}$\\
Now, for $\xi$ fixed, we  introduce the following time-Sobolev space~:
\begin{equation}\label{ForcingSpace}
Y^b_\xi=\{u\in\mathcal{S}'(\mathbb R^3);||u||_{Y^b_\xi}=:||\langle i\tau +\xi^2\rangle^b
 \hat u(\tau)||_{L^2_\tau(\mathbb R)}<\infty\}.
\end{equation}
In order to obtain certain estimates in $X^{b,s_1,s_2}$ for the following operator
\begin{equation}\label{ForcingOp2}
L:f\longmapsto\chi_{\mathbb R_{+}}(t) \psi(t)\int_0^t W(t-t')f(t')dt'
\end{equation}
we shall study in $Y^b_\xi$ the following linear operator~:
\begin{equation}\label{ForcingOp1}
K:f\longmapsto \psi(t)\int_0^t e^{-|t-t'|\xi^2}f(t')dt'
\end{equation}
\begin{proposition}\label{ProForcing1}
Let $\xi\in \mathbb R$ fixed and $f\in \mathcal S(\mathbb R^3)$, $0<\delta\leq 1/2$.
One considers the operator
\begin{equation}\label{ProForcing1Op}
t\longmapsto K_\xi(t)= \psi(t)\int_0^te^{-|t-t'|\xi^2} f(t')dt'.
\end{equation}
Then we have the following estimate~:
\begin{equation}\label{ProForcing1Es}
||K_\xi(t)||_{Y_\xi^{1/2}}\leq C\langle \xi\rangle^{-2\delta}||f||_{Y_\xi^{-1/2+\delta}}
\end{equation}\hfill$\Box$
\end{proposition}
$\bold{Proof.}$
A simple calculation in \cite{MoRi} gives,
\begin{equation}\label{ProForcMoRi}
K_\xi(t)=\psi(t)\int_{\mathbb R}\frac{e^{it\tau}-e^{-|t|\xi^2}}{i\tau+\xi^2}\hat f(\tau)d\tau.
\end{equation}
 We  can break up $K_\xi$ in $K_\xi=K_{1,0}+K_{1,\infty}+K_{2,0}+K_{2,\infty}$, where
\begin{equation}
\nonumber K_{1,0}=:\psi(t)\int_{|\tau|\leq 1}\frac{e^{it\tau}-1}{i\tau+\xi^2}\hat f(\tau)d\tau,
\hspace{0.2 cm}
K_{1,\infty}=\psi(t)\int_{|\tau|\geq 1}\frac{e^{it\tau}}{i\tau+\xi^2}\hat f(\tau)d\tau,
\end{equation}
\begin{equation}
\nonumber K_{2,0}=\psi(t)\int_{|\tau|\leq 1}\frac{1-e^{-|t|\xi^2}}{i\tau+\xi^2}\hat f(\tau)d\tau,
\hspace{0.2 cm}
K_{2,\infty}=\psi(t)\int_{|\tau|\geq 1}\frac{e^{-|t|\xi^2}}{i\tau+\xi^2}\hat f(\tau)d\tau.
\end{equation}
\underline{Contribution of $K_{1,0}$.} In this case, while using the  asymptotic expansion,
we have,
\begin{equation}\label{ProofForTay}
K_{1,0}=\psi(t)\sum_{n\geq 1}\int_{|\tau|\leq 1}\frac{(it\tau)^n}{i\tau+\xi^2}\hat f(\tau)d\tau
\end{equation}
it results that
\begin{eqnarray}\label{ProofForRes}
\nonumber||\langle i\tau +\xi^2\rangle^{1/2} \mathcal{F}_t( K_{1,0})||_{L^2_\tau(\mathbb R)}
\nonumber&\leq& \sum_{n\geq1}||\langle i\tau +\xi^2\rangle^{1/2}
\mathcal{F}_t(\frac{\psi(t)t^n}{n!})||_{L^2_\tau(\mathbb R)}
\\
\nonumber&&\times\int_{|\tau|\leq1}\frac{|i\tau|^n}{|i\tau+\xi^2|}|\hat f(\tau)|d\tau \\
\nonumber&\leq&\big{(}\sum_{n\geq1}||\frac{\psi(t)t^n}{n!}||_{H^{1/2}_\tau} + |\xi|
||\frac{t^n\psi(t)}{n!}||_{L^2_\tau}\big{)}
\\
&&\times\int_{|\tau|\leq1}\frac{|\tau|^n}{|i\tau+\xi^2|}|\hat f(\tau)|d\tau,
\end{eqnarray}
 Using the inequality  $||t^n\psi(t)||_{H^b}\leq C n$ for $b\in\{0,1/2\}$, $n\geq 1$, together with
  the Cauchy-Schwarz inequality, we obtain
\begin{eqnarray}\label{ProofForCh1}
\nonumber||K_{1,0}||_{Y^{1/2}_\xi}
\nonumber&\leq& C(1+|\xi|)(\int_{|\tau|\leq1}\frac{|\hat f(\tau)|^2}
{\langle i\tau +\xi^2\rangle}d\tau)^{1/2}
(\int_{|\tau|\leq1}\frac{|\tau|^2\langle i\tau +\xi^2\rangle}{|i\tau+\xi^2|^2}d\tau)^{1/2}
\\
\nonumber&\leq&C\langle\xi\rangle(\int_{|\tau|\leq1}\frac{|\hat f(\tau)|^2}{\langle i\tau +
\xi^2\rangle}d\tau)^{1/2}\langle\xi\rangle^{-1}
\\
&\leq& C(\int_{|\tau|\leq1}\frac{|\hat f(\tau)|^2}{\langle i\tau +\xi^2\rangle}d\tau)^{1/2}.
\end{eqnarray}
Finally, since for $0\leq\delta< 1/2$ we have
 $\langle i\tau +\xi^2\rangle\geq\langle i\tau +\xi^2\rangle^{1-2\delta}\langle\xi\rangle^{4\delta}$,
it results that
\begin{eqnarray}\label{ProofForFin}
\nonumber||K_{1,0}||_{Y^{1/2}_\xi}\nonumber&\leq&C\langle\xi\rangle^{-2\delta}(\int_{|\tau|\leq1}
\frac{|\hat f(\tau)|^2}{\langle i\tau +\xi^2\rangle^{1-2\delta}}d\tau)^{1/2}
\\
&\leq&C\langle\xi\rangle^{-2\delta}||f||_{Y^{-1/2+\delta}_\xi}.
\end{eqnarray}
\underline{Contribution of $K_{2,\infty}$.}
 We note that,
\begin{eqnarray}\label{ProofForNot2}
\nonumber||\langle i\tau +\xi^2\rangle^{1/2} \mathcal{F}_t( K_{2,\infty})||_{L^2_\tau(\mathbb R)}
\nonumber&\leq& ||\langle i\tau +\xi^2\rangle\mathcal{F}_t(\psi(t)e^{-\xi^2|t|})
||_{L^2_\tau(\mathbb R)}\\
&&\times (\int_{|\tau|\geq1}\frac{|\hat f(\tau)|}{\langle i\tau +\xi^2\rangle}d\tau),
\end{eqnarray}
 Using the  inequality (\ref{FreeMoRi}), get now that,
\begin{equation}\label{ProofForFree2}
 ||\langle i\tau +\xi^2\rangle^{1/2}\mathcal{F}_t(\psi(t)e^{-\xi^2|t|})||_{L^2_\tau(\mathbb R)}\leq C,
 \end{equation}
therefore, by the Cauchy Schwarz inequality we obtain,
\begin{eqnarray}\label{ProofForRes2}
\nonumber||K_{2,\infty}||_{Y_\xi^{1/2}}&\leq C& (\int_{|\tau|\geq1}\frac{|\hat f(\tau)|}
{\langle i\tau +\xi^2\rangle}d\tau)\\
&\leq& C (\int_{\mathbb R}\frac{|\hat f(\tau)|^2}{{\langle i\tau +\xi^2\rangle}^{1-2\delta}}
d\tau)^{1/2} (\int_{|\tau|\geq1}\frac{\langle i\tau +\xi^2\rangle^{-2\delta}}{\langle i\tau +
\xi^2\rangle}d\tau)^{1/2}
\end{eqnarray}
 For $| \xi|\geq 1$,  the following change of variable $\tau\longmapsto r\xi^2$, give

\begin{eqnarray}\label{ProofForCV2}
\nonumber||K_{2,\infty}||_{Y_\xi^{1/2}}&\leq C& (\int_{\mathbb R}\frac{|\hat f(\tau)|^2}
{{\langle i\tau +\xi^2\rangle}^{1-2\delta}}d\tau)^{1/2}\langle \xi\rangle^{-2\delta}
(\int_{\mathbb R}\frac{1}{\langle r\rangle^{1+2\delta}}dr)^{1/2}\\
&\leq&C \langle \xi\rangle ^{-2\delta}||f||_{Y_\xi^{-1/2+\delta}}
\end{eqnarray}
since we have  $\displaystyle\int_{\mathbb R}\frac{1}{\langle r\rangle^{1+2\delta}}dr < \infty$.\\
In the other case when $|\xi|\leq1$ it follows $\langle \xi\rangle ^{-2\delta}\sim 1$ . Therefore
\begin{eqnarray}\label{ProofForCV21}
\nonumber||K_\gamma^{2,\infty}||_{Y_{\xi}^{1/2}}&\lesssim&\langle \xi\rangle ^{-2\delta}  (\int_{\mathbb R}\frac{|\hat f(\tau)|^2}{{\langle i\tau +\xi^2\rangle}^{1-2\delta}}d\tau)^{1/2}
(\int_{\mathbb R}\frac{1}{\langle \tau\rangle^{1+2\delta}}d)^{1/2}\\
&\lesssim&\langle \xi\rangle ^{-2\delta}||f||_{Y_{\gamma,\xi}^{-1/2+\delta}}
\end{eqnarray}

\underline{Contribution of $K_{2,0}$}. We note that,
\begin{equation}\label{ProofForNot3}
||K_{2,0}||_{Y_\xi^{1/2}}\leq ||\langle i\tau +\xi^2\rangle ^{1/2}\mathcal{F}_t(\psi(t)
(1-e^{-\xi^2|t|}))(\tau)||_{L^2_\tau}\int_{|\tau|\leq1}\frac{|\hat f(\tau)|}{|i\tau +\xi^2|}d\tau.
\end{equation}
Case 1: $|\xi|\geq1$. Using the inequality (\ref{FreeMoRi}) in the proof of  Propositon \ref{ProFree}
we obtain,
\begin{eqnarray}\label{ProofForMoRi31}
\nonumber I &=:& ||\langle i\tau +\xi^2\rangle ^{1/2}\mathcal{F}_t(\psi(t)
(1-e^{-\xi^2|t|}))(\tau)||_{L^2_\tau}\\
\nonumber&\leq&||\psi(t)(1-e^{-\xi^2|t|})(\tau)||_{H^{1/2}_\tau}
+ \langle\xi \rangle||\mathcal{F}_t(\psi(t)(1-e^{-\xi^2|t|}))(\tau)||_{L^2_\tau}\\
&\leq& C(1+\langle \xi \rangle)\leq C\langle \xi \rangle,
\end{eqnarray}
therefore,
\begin{equation}\label{ProofForTher31}
||K_{2,0}||_{Y_\xi^{1/2}}\leq C\langle \xi \rangle\int_{|\tau|\leq1}
\frac{|\hat f(\tau)|}{|i\tau +\xi^2|}d\tau,
\end{equation}
now, we apply  the  Cauchy Schwarz inequality to obtain,
\begin{eqnarray}\label{ProofForCS31}
\nonumber||K_{2,0}||_{Y_\xi^{1/2}}&\leq& C\langle \xi \rangle(\int_{|\tau|\leq1}
\frac{|\hat f(\tau)|^2}{\langle i\tau +\xi^2\rangle}d\tau)^{1/2}(\int_{|\tau|\leq1}
\frac{\langle i\tau+\xi^2\rangle}{|i\tau +\xi^2|^2}d\tau)^{1/2}\\
\nonumber&\leq&C\langle \xi \rangle|\xi|^{-1}(\int_{|\tau|\leq1}
\frac{|\hat f(\tau)|^2}{\langle i\tau +\xi^2\rangle}d\tau)^{1/2}\\
\nonumber&\leq&C(\int_{|\tau|\leq1}\frac{|\hat f(\tau)|^2}{\langle i\tau +\xi^2\rangle}d\tau)^{1/2}\\
\nonumber&\leq&C\langle \xi \rangle^{-2\delta}(\int_{|\tau|\leq1}
\frac{|\hat f(\tau)|^2}{\langle i\tau +\xi^2\rangle^{1-2\delta}}d\tau)^{1/2}\\
&\leq& C\langle\xi\rangle^{-2\delta}||f||_{Y_\xi^{-1/2+\delta}}
\end{eqnarray}
Case 2~: $|\xi|\leq 1$. In this case we note that,
\begin{eqnarray}\label{ProofForNot32}
\nonumber I&\leq &||\psi(t)(1-e^{-\xi^2|t|})(\tau)||_{H^{1/2}_\tau}
+ \langle\xi \rangle||\mathcal{F}_t(\psi(t)(1-e^{-\xi^2|t|}))(\tau)||_{L^2_\tau}\\
\nonumber&\leq&||\psi(t)(1-e^{-\xi^2|t|})(\tau)||_{H^{1/2}_\tau}
\leq\sum_{n\geq1}||\frac{|t|^n\psi(t)|\xi|^{2n}}{n!}||_{H^{1/2}_\tau}\\
&\leq&C\sum_{n\geq1}\frac{|\xi|^{2n}}{n!}|||t|^{n}\psi(t)||_{H^{1/2}_\tau}
\leq C|\xi|\sum_{n\geq1}\frac{1}{(n-1)!}\leq C|\xi|.
\end{eqnarray}
Substituting this inequality in (\ref{ProofForNot3}), then  as in Case 1, we obtain using  Cauchy
Schwarz inequality ,
\begin{eqnarray}\label{ProofForSub3}
\nonumber ||K_{2,0}||_{Y^{1/2}_\xi}&\leq& C|\xi||\xi|^{-1}(\int_{|\tau|\leq1}
\frac{|\hat f(\tau)|^2}{\langle i\tau +\xi^2\rangle}d\tau)^{1/2}\\
&\leq& C\langle\xi\rangle^{-2\delta}||f||_{Y_\xi^{-1/2+\delta}}
\end{eqnarray}
\underline{Contribution of $K_{1,\infty}$}. By the identity $\mathcal F(u\star v)=\hat u\hat v$
and the triangle inequality $\langle i\tau + \xi^2\rangle \leq \langle \tau_1 \rangle +|i(\tau-\tau_1)+\xi^2|$, we see that
\begin{eqnarray}\label{ProofForTr4}
\nonumber||K_{1,\infty}||_{Y^{1/2}_\xi}&=&||\langle i\tau +\xi^2\rangle^{1/2}
|\hat\psi(\tau_1)\star(\frac{\hat f(\tau_1)}{|i\tau_1 +\xi^2|}\chi_{\{|\tau_1|\geq 1\}})|(\tau)||_{L^2_\tau}\\
\nonumber&\leq&||\hspace{0,1cm}|\langle \tau_1 \rangle^{1/2}\hat\psi(\tau_1)|\star(\frac{|\hat f(\tau_1)|}
{|i\tau_1 +\xi^2|}\chi_{\{|\tau_1|\geq 1\}})(\tau)||_{L^2_\tau}\\
&&+||\hspace{0,1cm}|\hat\psi(\tau_1)|\star(\frac{|\hat f(\tau_1)|}{|i\tau_1 +\xi^2|^{1/2}}\chi_{\{|\tau_1|\geq 1\}})
(\tau)||_{L^2_\tau}.
\end{eqnarray}
Due to the convolution inequality $||u\star v||_{L^2_\tau}\lesssim ||u||_{L^1_\tau}||v||_{L^2_\tau}$,
 we obtain,
\begin{eqnarray}\label{ProofForLC4}
\nonumber||K_{1,\infty}||_{Y^{1/2}_\xi}&\leq&||\langle \tau \rangle\hat\psi(t)||_{L^1_{\tau}}||
\frac{|\hat f(\tau)|}{|i\tau +\xi^2|}\chi_{\{|\tau|\geq 1\}}||_{L^2_\tau}\\
\nonumber&&+||\psi(t)||_{L^1_\tau}||\frac{|\hat f(\tau)|}{|i\tau +\xi^2|^{1/2}}
\chi_{\{|\tau|\geq 1\}}||_{L^2_\tau}\\
\nonumber&\leq& C||\frac{|\hat f(\tau)|}{\langle i\tau +\xi^2\rangle^{1/2}}
\chi_{\{|\tau|\geq 1\}}||_{L^2_\tau}\\
\nonumber&\leq& C \langle \xi \rangle^{-2\delta}||\langle i\tau +\xi^2\rangle^{-1/2+\delta}
\hat f(\tau)||_{L^2_\tau}\\
&\leq&C \langle \xi \rangle^{-2\delta}||f||_{Y_\xi^{-1/2+\delta}}.
\end{eqnarray}
This completes the proof of the Proposition.$\hfill\Box$\\\\
Now, by  use of Proposition \ref{ProForcing1}, we prove some smoothing properties in the Bourgain
 spaces for the operator $L$  defined by (\ref{ForcingOp2}).

\begin{proposition}\label{PrOp2}
Let $0<\delta\leq 1/2$, $s_1$ and $s_2\in \mathbb R$, there exists $C=C(\delta)>0$ such that,
 for all $f\in X^{-1/2+\delta, s_1-2\delta,s_2}$, we have
\begin{equation}
M=:||\chi_{\mathbb R_{+}}(t) \psi(t)\int_0^t W(t-t')f(t')dt'||_{X^{1/2,s_1,s_2}}\leq
 C_\delta||f||_{X^{-1/2 +\delta,s_1-2\delta,s_2}}
\end{equation}$\hfill\Box$
\end{proposition}
$\bold{Proof}$. By definition of $X^{1/2,s_1,s_2}$, we see that $M$ is equal to
\begin{equation}\label{ProOp2D}
||\langle \xi\rangle^{s_1}\langle \eta\rangle^{s_2}\langle i(\tau-P(\nu))+\xi^2\rangle^{1/2}
\mathcal{F}_{t,x}(\chi_{\mathbb R_{+}}(t) \psi(t)\int_0^t W(t-t')f(t')dt')
(\tau,\nu)||_{L^2_{\tau,\nu}},
\end{equation}
we note that,
$$\hspace{-4cm}\mathcal{F}_{t,x}(\chi_{\mathbb R_{+}}(t) \psi(t)\int_0^t W(t-t')f(t')dt)(\tau)$$
\begin{eqnarray*}
&=&\mathcal{F}_{t}(\chi_{\mathbb R_{+}}(t) \psi(t)\int_0^t e^{-|t-t'|\xi^2}e^{it'P(\nu)(t-t')}\mathcal{F}_\nu(f)(t',\nu)dt')(\tau)\\
&=&\mathcal{F}_{t}(\chi_{\mathbb R_{+}}(t) \psi(t)\int_0^t e^{-|t-t'|\xi^2}e^{-it'P(\nu)t'}\mathcal{F}_\nu(U(t)f(t',\nu))dt')(\tau)\\
&=&\mathcal{F}_{t}(\chi_{\mathbb R_{+}}(t) \psi(t)\int_0^t e^{-|t-t'|\xi^2}e^{-it'P(\nu)t'}\mathcal{F}_\nu(f)(t',\nu)dt')\big(\tau-P(\nu)\big),
\end{eqnarray*}
then by the change of variable $\tau\longmapsto \tau-P(\nu)$, we obtain
$$
 M=||\langle \xi\rangle^{s_1}\langle \eta\rangle^{s_2}
\mathcal{F}_{t}(\chi_{\mathbb R_{+}}(t) \psi(t)\int_0^t e^{-|t-t'|\xi^2}e^{-it'P(\nu)t'}\mathcal{F}_\nu(f)(t',\nu)dt')||_{L^2_{\nu}(Y^{1/2}_\xi)}$$
\begin{equation}\label{ProOp2Ch}
=||\langle \xi\rangle^{s_1}\langle \eta\rangle^{s_2}
\mathcal{F}_{t}(\chi_{\mathbb R_{+}}(t) \psi(t)\int_0^t e^{-|t-t'|\xi^2}\mathcal{F}_\nu(U(-t')f)(t',\nu)dt')||_{L^2_{\nu}(Y^{1/2}_\xi)},
\end{equation}
now, let us set $w(\tau,\nu)=\mathcal{F}_\nu(U(-t')f)(\tau,\nu)$. To apply the Proposition
\ref{ProForcing1}, we  need firstly to suppose that $f \in \mathcal S(\mathbb R^3)$. It is clear that
 $w\in \mathcal S(\mathbb R^3)$ and we take
\begin{equation}
 K_\xi:f\longmapsto \psi(t)\int_0^t e^{-|t-t'|\xi^2}w(t')dt,
\end{equation}
therefore,
\begin{eqnarray}\label{ProOp2K}
\nonumber
 M&=&||\langle \xi\rangle^{s_1}\langle \eta\rangle^{s_2}
\mathcal{F}_{t}(\chi_{\mathbb R_{+}}(t) K_\xi(t))||_{L^2_{\nu}(Y^{1/2}_\xi)}\\
&\leq&||\langle \xi\rangle^{s_1}\langle \eta\rangle^{s_2}
||\chi_{\mathbb R_{+}}(t) K_\xi(t)||_{H^{1/2}_\tau}||_{L^2_{\nu}}\\
&&+||\langle \xi\rangle^{s_1+1}\langle \eta\rangle^{s_2}
||\chi_{\mathbb R_{+}}(t) K_\xi(t)||_{L^2_\tau}||_{L^2_{\nu}}.
\end{eqnarray}
Since $K_\xi(0)=0$, then we have  $||\chi_{\mathbb R_{+}}(t) K_\xi(t)||_{H^{1}_\tau}\leq|| K_\xi(t)||_{H^{1}_\tau}$ and  $||\chi_{\mathbb R_{+}}(t) K_\xi(t)||_{L^2_\tau} \leq|| K_\xi(t)||_{L^2_\tau}$
, by noting that $ ||h||_{L^2}=||h||_{H^0}$, it results  by interpolation between  $H^0$ and $ H^1$
that
\begin{equation}\label{ProOp2Int}
 ||\chi_{\mathbb R_{+}}(t) K_\xi(t)||_{H^{1/2}_\tau}\leq|| K_\xi(t)||_{H^{1/2}_\tau}.
\end{equation}
It results  that
\begin{eqnarray}\label{ProOp2RM}
\nonumber M&\leq&||\langle \xi\rangle^{s_1}\langle \eta\rangle^{s_2}
|| K_\xi(t)||_{H^{1/2}_\tau}||_{L^2_{\nu}}
+||\langle \xi\rangle^{s_1+1}\langle \eta\rangle^{s_2}
|| K_\xi(t)||_{L^2_\tau}||_{L^2_{\nu}}\\
 &\leq& C||\langle \xi\rangle^{s_1}\langle \eta\rangle^{s_2}
\mathcal{F}_{t}( K_\xi(t))||_{L^2_{\nu}(Y^{1/2}_\xi)},
\end{eqnarray}

now, we can apply the Proposition \ref{ProForcing1} to obtain,
 \begin{eqnarray}\label{ProOp2Forc}
\nonumber M&\leq&C||\langle \xi\rangle^{s_1}\langle \eta\rangle^{s_2}
||w||_{Y^{-1/2+\delta}_\xi}\langle \xi\rangle^{-2\delta}||_{L^2_{\nu}}\\
\nonumber&=&C||\langle \xi\rangle^{s_1-2\delta}\langle \eta\rangle^{s_2}
||\langle i\tau +\xi^2\rangle ^{-1/2+\delta}\mathcal{F}_{t}( w(t))||_{L^2_\tau}||_{L^2_{\nu}}\\
\nonumber&=&C||\langle \xi\rangle^{s_1-2\delta}\langle \eta\rangle^{s_2}
\langle i\tau +\xi^2\rangle ^{-1/2+\delta}\mathcal{F}_{t,x}( U(-t)f(\tau,\nu))||_{L^2_{\tau,\nu}}\\
&=&C||\langle \xi\rangle^{s_1-2\delta}\langle \eta\rangle^{s_2}
\langle i\tau +\xi^2\rangle ^{-1/2+\delta} \hat f(\tau + P(\nu),\nu)||_{L^2_{\tau,\nu}},
\end{eqnarray}
finally, by the change of variable $ \tau \longmapsto \tau -P(\nu)$ we  can deduce that
$M\leq C||f||_{X^{-1/2+\delta,s_1-2\delta,s_2}}$, this for any $ f\in
\mathcal S(\mathbb R^3).$ The result for $f\in X^{-1/2+\delta,s_1-2\delta,s_2}$ follows by density.
$\hfill\Box$
\begin{proposition}\label{PrForReg}
Let $s_1,$ $s_2\in \mathbb R$ and $0< \delta\leq 1/2$. For all $f\in X^{-1/2+\delta, s_1-2\delta,s2}$ we have
\begin{equation}\label{PrForReg1}
L:t\longmapsto \int_0^tW(t-t')f(t')dt'\in C(\mathbb R_+, H^{s_1,s_2})
\end{equation}
moreover, if $(f_n)$ is a sequence with $f_n\longrightarrow 0 \text{ in } X^{-1/2+\delta, s_1-2\delta,s2}$
 if $n\to\infty$,then
\begin{equation}\label{PrForReg2}
|| \int_0^tW(t-t')f_n(t')dt'||_{L^{\infty}(\mathbb R_+, H^{s_1,s_2})}\longrightarrow 0
\end{equation}
$\hfill\Box$
\end{proposition}

$\bold{Proof}$. By  Fubini theorem, and by the definition of $W(.)$ we have,
\begin{eqnarray}
\nonumber L(t) &=&  \int_0^t\mathcal F^{-1}_\nu (e^{-|t-t'|\xi^2}e^{i(t-t')P(\nu)}
\mathcal F_\nu(f(t')))dt'\\
\nonumber &=&\int_{\mathbb R^2}e^{i\langle \nu,(x,y)\rangle}e^{itP(\nu)}
 \int_0^t e^{-|t-t'|\xi^2}\mathcal F_{(x,y)}(U(-t')f(t'))(\nu)dt'\\
&=&U(t)\mathcal F^{-1}_{(x,y)}\Big[ \int_0^t e^{-|t-t'|\xi^2}\mathcal F_{(x,y)}(g(t',.))(\nu)dt'\Big]
\end{eqnarray}
where $g(t,\nu)=(U(-t)f(t))(\nu)$. As noticed in \cite{Ge95} since $U(\cdot)$ is a strongly continuous unitary group in$L^2(\mathbb R^2)$, it is enough to prove that
$t\longmapsto U(-t)L(t)\in C(\mathbb R_{+},H^{s_1,s_2})$, then it is equivalent to show that,
\begin{equation}\label{PrForReg3}
F:t\longmapsto \langle \xi\rangle^{s_1}\langle \eta\rangle^{s_2} \int_0^t e^{-|t-t'|\xi^2}
\mathcal F_{(x,y)}(g(t',.))(\nu)dt'
\end{equation}
 is continuous from $ \mathbb R_{+}\longmapsto L^2(\mathbb R^2)$, for
$f\in X^{-1/2+\delta, s_1 -2\delta,s_2}$, $0<\delta\leq 1/2$. Note that by the Fubini theorem
we have,
\begin{eqnarray}\label{PrForReg4}
\nonumber F(t)&=& \langle \xi\rangle^{s_1}\langle \eta\rangle^{s_2}e^{-t|\xi|^2} \int_0^t e^{t'|\xi|^2}\mathcal F_{(x,y)}(g(t',.))(\nu)dt'\\
\nonumber&=& \langle \xi\rangle^{s_1}\langle \eta\rangle^{s_2}e^{-t|\xi|^2}\int_{\mathbb R}\hat g(\tau, \nu) \int_0^t e^{(i\tau+|\xi|^2)t'}dt'd\tau\\
&=& \langle \xi\rangle^{s_1}\langle \eta\rangle^{s_2}\int_{\mathbb R}\hat g(\tau, \nu) \frac{e^{it\tau}-e^{-|\xi|^2t}}{i\tau+\xi^2} d\tau .
\end{eqnarray}
one fixes $t_1$, $t_2\in \mathbb R_{+}$, then,
\begin{eqnarray}\label{PrForReg5}
\nonumber F(t_1)-F(t_2)&=& \langle \xi\rangle^{s_1}\langle \eta\rangle^{s_2}
\int_{\mathbb R} \frac{\hat g(\tau, \nu)}{i\tau+\xi^2}\Big[(e^{it_1\tau}-e^{it_2\tau})\\
\nonumber&&\qquad -(e^{-|\xi|^2t_1}-e^{-|\xi|^2t_2})\Big] d\tau \\
&=:& \langle \xi\rangle^{s_1}\langle \eta\rangle^{s_2}
\int_{\mathbb R} J_{t_1,t_2}(\tau) d\tau.
\end{eqnarray}
We  deal first with the case   $|\xi|\geq1$.  Using  Cauchy Schwarz inequality we obtain,

\begin{eqnarray}\label{PrForReg6}
\nonumber| F(t_1)-F(t_2)|&\leq& 4\langle \xi\rangle^{s_1}\langle \eta\rangle^{s_2}
\int_{\mathbb R} \frac{|\hat g(\tau, \nu)|}{|i\tau+\xi^2|}d\tau\\
&\leq& 4\langle \xi\rangle^{s_1}\langle \eta\rangle^{s_2}
\Big(\int_{\mathbb R} \frac{|\hat g(\tau, \nu)|^2}{\langle i\tau+\xi^2\rangle^{1-2\delta}}d\tau
\Big)^{1/2}
\Big(\int_{\mathbb R} \frac{\langle i\tau+\xi^2\rangle^{1-2\delta}}{|i\tau+\xi^2|^2}d\tau\Big)^{1/2}.
\end{eqnarray}
Since in this case we have $|i\tau+\xi^2|\sim\langle i\tau +\xi^2\rangle$, by the change of variable
:~$\tau\longmapsto r\xi^2$ it results,
\begin{eqnarray}\label{PrForReg7}
\nonumber| F(t_1)-F(t_2)|&\leq& C\langle \xi\rangle^{s_1}\langle \eta\rangle^{s_2}
\Big(\int_{\mathbb R} \frac{|\hat g(\tau, \nu)|^2}{\langle i\tau+\xi^2\rangle^{1-2\delta}}d\tau \Big)^{1/2}
|\xi|^{-2\delta}\Big(\int_{\mathbb R} \frac{dr}{\langle
 r\rangle^{1+2\delta}}d\tau\Big)^{1/2}\\
&\leq& C\langle \xi\rangle^{s_1-2\delta}\langle \eta\rangle^{s_2}||\langle i\tau +\xi^2\rangle^{-1/2+\delta}\hat g(.,\nu)||_{L^2_\tau}.
\end{eqnarray}
In the other case when $|\xi|\leq1$,  we assume that $|t_1-t_2|$ is small enough. We can write
$ |F(t_1)-F(t_2)|\leq I_1+I_2 $, where
\begin{equation}
I_1=:\langle \xi\rangle^{s_1}\langle \eta\rangle^{s_2}
\Big|\int_{\mathbb R} \frac{\hat g(\tau, \nu)}{i\tau+\xi^2}\big[e^{it_1\tau}-e^{it_2\tau}\big]
 d\tau \Big|,
\end{equation}
\begin{equation}
I_2=:\langle \xi\rangle^{s_1}\langle \eta\rangle^{s_2}
\Big|\int_{\mathbb R} \frac{\hat g(\tau, \nu)}{i\tau+\xi^2}\big[e^{-\xi^2t_1}-e^{-\xi^2t_2}\big]
 d\tau \Big |.
\end{equation}
Now, we start to estimate  $I_1$. By Cauchy Schwarz inequality, we see that
\begin{eqnarray}\label{PrForReg8}
\nonumber I_1&=:&\langle \xi\rangle^{s_1}\langle \eta\rangle^{s_2}
[|t_1-t_2|\int_{|\tau|\leq1} |\tau|\bigg|\frac{\hat g(\tau, \nu)}{i\tau+\xi^2}\bigg|d\tau +
2\int_{|\tau|\geq1} \bigg|\frac{\hat g(\tau, \nu)}{i\tau+\xi^2}\bigg|d\tau]\\
\nonumber&\leq&C\langle \xi\rangle^{s_1}\langle \eta\rangle^{s_2}||\langle i\tau +\xi^2\rangle^{-1/2+\delta}\hat g(.,\nu)||_{L^2_\tau}\bigg[\Big(\int_{|\tau|\leq1} \frac{|\tau|^2\langle i\tau+\xi^2\rangle^{1-2\delta}}{|i\tau+\xi^2|^2}d\tau\Big)^{1/2}
\\
&& \hspace{3 cm}+
\Big(\int_{|\tau|\geq1} \frac{|i\tau+\xi^2|^{1-2\delta}}{|i\tau+\xi^2|}d\tau\Big)^{1/2}\bigg].
\end{eqnarray}
Since for  $|\xi|\leq1$, we have  $\langle i\tau +\xi^2\rangle^{1-2\delta}\sim \langle\tau
\rangle^{1-2\delta}$. Using this approximation in the first integral of (\ref{PrForReg8}) together
 with   the change of variable :~$\tau\longmapsto r\xi^2$ in  the second one, we obtain
\begin{eqnarray}\label{PrForReg9}
\nonumber I_1&\leq&C\langle \xi\rangle^{s_1}\langle \eta\rangle^{s_2}||\langle i\tau +\xi^2\rangle^{-1/2+\delta}\hat g(.,\nu)||_{L^2_\tau}\Big((\int_{|\tau|\leq1} \langle \tau\rangle^{1-2\delta}d\tau)^{1/2}
\\
\nonumber&& \hspace{3 cm}+
(\int_{\mathbb R} |\frac{1}{\langle r\rangle ^{1+2\delta}}dr)^{1/2})\Big)\\
&\leq&C\langle \xi\rangle^{s_1}\langle \eta\rangle^{s_2}||\langle i\tau +\xi^2\rangle^{-1/2+\delta}\hat g(.,\nu)||_{L^2_\tau}
\end{eqnarray}
Note that in this case  $\langle \xi \rangle^{-2\delta}\sim 1$, therefore
\begin{equation}\label{PrForReg10}
I_1\leq C\langle \xi\rangle^{s_1-2\delta}\langle \eta\rangle^{s_2}||\langle i\tau +\xi^2\rangle^{-1/2+\delta}\hat g(.,\nu)||_{L^2_\tau}
\end{equation}
Now, we pass to estimate  $I_2$. By Cauchy Schwarz inequality it results that,
\begin{equation}\label{PrForReg11}
I_2\leq C\langle \xi\rangle^{s_1}\langle \eta\rangle^{s_2}|\xi|^2||\langle i\tau +
\xi^2\rangle^{-1/2+\delta}\hat g(.,\nu)||_{L^2_\tau}(\int_{\mathbb R}\frac{\langle i\tau+
\xi^2\rangle ^{1-2\delta}}{|i\tau+ \xi^2|^2}d\tau)^{1/2}.
\end{equation}
Since we have,
\begin{eqnarray}\label{PrForReg12}
\nonumber\int_{\mathbb R}\frac{\langle i\tau+\xi^2\rangle ^{1-2\delta}}{|i\tau+ \xi^2|^2}d\tau&=&
\int_{|\tau|\leq1}\frac{\langle i\tau+\xi^2\rangle ^{1-2\delta}}{|i\tau+ \xi^2|^2}d\tau +
\int_{|\tau|\geq1}\frac{\langle i\tau+\xi^2\rangle ^{1-2\delta}}{|i\tau+ \xi^2|^2}d\tau\\
&\leq& C(|\xi|^{-4}+|\xi|^{-4\delta})
\end{eqnarray}
then,
\begin{eqnarray}\label{PrForReg13}
\nonumber I_2 &\leq& C\langle \xi\rangle^{s_1}\langle \eta\rangle^{s_2}|\xi|^2(|\xi|^{-2}+|\xi|^{-2\delta})||\langle i\tau +\xi^2\rangle^{-1/2+\delta}\hat g(.,\nu)||_{L^2_\tau}\\
 &\leq& C\langle \xi\rangle^{s_1-2\delta}\langle \eta\rangle^{s_2}||\langle i\tau +\xi^2\rangle^{-1/2+\delta}\hat g(.,\nu)||_{L^2_\tau}
\end{eqnarray}
finally, gathering (\ref{PrForReg7}), (\ref{PrForReg10}),(\ref{PrForReg13}), one infers that
\begin{eqnarray}\label{PrForReg14}
\nonumber||F(t_1)-F(t_2)||_{L^2_\tau}&\leq& C ||\langle \xi\rangle^{s_1-2\delta}\langle \eta\rangle^{s_2}\langle i\tau +\xi^2\rangle^{-1/2+\delta}\hat g(.,\nu)||_{L^2_{\tau,\nu}}\\
&=&C||f||_{X^{-1/2+\delta,s_1-2\delta,s_2}}.
\end{eqnarray}
It is clear that the integrand in (\ref{PrForReg5}) tends to $0$ pointwise in $(\tau, \nu)$ as soon as $|t_1-t_2|\to 0$ and is bounded uniformly in $|t_1-t_2|$ by  the right member of (\ref{PrForReg14}). the result follows then from Lebesgue dominated convergence theorem.\\
To show (\ref{PrForReg2})it suffices to notice that one has
\begin{equation}\label{PrForReg15}
\nonumber\sup_{t\in\mathbb R_{+}}||F_n(t)||_{L^2(\mathbb R^2)}\leq C||f_n||_{X^{-1/2+\delta,s_1-2\delta,s_2}}
\end{equation}
where $F_n$ is defined as $F$ with $ g_n(.)=\mathcal F_\nu(U(-t)f_n(t))$ instead of g. This completes the proof. $\hfill{\Box}$

\section{Strichartz  and multilinear estimates for the KP-equation}\label{EstStr}
The goal in this section is to prepare certain  Strichartz and multilinear  estimates  by using
result derived  by Molinet-Ribaut in \cite{MoRi} and  Saut in \cite{Sau93}
This type of estimates is necessary to treat in the next section the
nonlinear term $\partial(u^2)$ in $X^{b,s_1,s_2}$. The following lemma is prepared by
Molinet-Ribaud in \cite{MoRi}.
\begin{lemma}\label{StrMoRi}
Let $v\in L^2(\mathbb R^2)$ whith $ supp v \subset \{(t,x,y): |t|\leq T\}$ and let $\epsilon>0$ ,
$\delta(r)=1-2/r$.
Then for all $(r,\beta,\theta)$ with
\begin{equation}\label{StrMoRi1}
2\leq r<\infty,\hspace {0,2 cm} 0\leq \beta \leq 1/2, \hspace{0,2 cm} 0\leq\delta (r)\leq \frac{\theta}{1-\beta/3}
\end{equation}
there exists $\mu=\mu(\epsilon)>0$ such that
\begin{equation}\label{StrMoRi2}
||\mathcal F^{-1}_{t,x}( |\xi|^{-\frac{\beta \delta(r)}{2}}\langle\tau- P(\nu)\rangle^{\frac{-\theta}{2} (1+\epsilon)}|\hat v(\tau, \nu)|)||_{L^{q,r}_{t,x}}\leq CT^\mu ||v||_{L^2(\mathbb R^3)}
\end{equation}
where $q$ is defined by
\begin{equation}\label{StrMoRi3}
2/q=(1-\beta/3)\delta(r) +(1-\theta)
\end{equation} $\hfill{\Box}$
\end{lemma}
Now, we will use  Lemma \ref{StrMoRi} to derive  a first  multilinear estimate.
\begin{lemma}\label{StrBass}
let $u$, $v$ with compact support in $\{(x,y,t):|t|\leq T\}$. For $b>0$ small enough, there exists $\mu >0$ such that
\begin{eqnarray}\label{StrBass1}
\nonumber I&=:& \int_{\mathbb R^6}\frac{|\hat u(\tau_1,\nu_1)||\hat v(\tau-\tau_1,\nu-\nu_1)||\hat w(\tau,\nu)|}
{\langle \sigma \rangle^{1/2-b}|\xi_1|^{b/4}\langle\sigma_1\rangle^b \langle\sigma_2\rangle^{1/2}}d\tau d\tau_1 d\nu d\nu_1\\
&&\hspace{2cm}\leq CT^\mu ||u||_{L^2_{t,x}}||v||_{L^2_{t,x}}||w||_{L^2_{t,x}}
\end{eqnarray}
where $\sigma$, $\sigma_1$ and $\sigma_2$ are defined by
\begin{equation}
\sigma=\tau-P(\nu),\hspace{0,2cm} \sigma_1=\tau_1-P(\nu_1), \hspace{0,2cm} \sigma_2=\tau-\tau_1 -P(\nu-\nu_1)
\end{equation}$\hfill{\Box}$
\end{lemma}
$\bold{Proof}.$
By the Plancherel Theorem we  see that
\begin{equation}\label{StrBass2}
I= \int_{\mathbb R^4}\mathcal F^{-1}_{t,x}(\frac{|\hat w(\tau,\nu)|}
{\langle \sigma \rangle^{1/2-b}})\mathcal F^{-1}_{t,x}(\frac{|\hat u(\tau,\nu)|}
{\langle \sigma \rangle^{b}|\xi|^{b/4}}\star\frac{|\hat v(\tau,\nu)|}
{\langle \sigma \rangle^{1/2}})(\tau,\nu) d\tau  d\nu,
\end{equation}
by using the fact that $\mathcal F^{-1}_{t,x}(h\star f)=\mathcal F^{-1}_{t,x}(h)\star\mathcal F^{-1}_{t,x}(f)$  then by applying Holder inequality in space and next in time we obtain that $I$ is bounded by the product of the three terms
\begin{equation}\label{StrBass3}
||\mathcal F^{-1}_{t,x}(\frac{|\hat w(\tau,\nu)|}{\langle \sigma \rangle^{1/2-b}})||_{L^{q_1,r_1}_{t,x}}
||\mathcal F^{-1}_{t,x}(\frac{|\hat u(\tau,\nu)|}{\langle \sigma \rangle^{b}|\xi|^{b/4}})||_{L^{q_2,r_2}_{t,x}}
||\mathcal F^{-1}_{t,x}(\frac{|\hat v(\tau,\nu)|}{\langle \sigma \rangle^{1/2}})||_{L^{q_3,r_3}_{t,x}}
\end{equation}
  where
\begin{equation}\label{StrBass4}
\sum_{i=1}^{3}1/q_i=1, \hspace{0,2cm} \sum_{i=1}^{3}1/r_i=1.
\end{equation}
Our goal now is to estimate the three terms of (\ref{StrBass3})  by using  Lemma \ref{StrMoRi}.
Let $b$ small enough we take first $\epsilon_1=\epsilon_2=\epsilon_3=\epsilon$, where
$\epsilon=\epsilon(b)$ will be a small paremetr. Also we choose
\begin{equation}\label{StrBass5}
\theta_1=\frac{1-2b}{1+\epsilon},\hspace{0.2 cm}\theta_2=\frac{2b}{1+\epsilon},\hspace{0.2 cm}
\theta_3=\frac{1}{1+\epsilon}
\end{equation}
we choose $\beta_1=\beta_3=0$. From \ref{StrMoRi3}, it remains to find $\beta_2$, $q_i$ and $r_i$ with
\begin{equation}\label{StrBass6}
\frac{2}{q_1}=\delta(r_1)+(1-\theta_1),\hspace{0.2 cm}\frac{2}{q_2}=\delta(r_2)+(1-\theta_2)-
\frac{\beta_2\delta(r_2)}{3},\hspace{0.2 cm}\frac{2}{q_3}=\delta(r_3)+(1-\theta_3),
\end{equation}
\begin{equation}\label{StrBass7}
\beta_2\delta(r_2)=\frac{b}{2}
\end{equation}
such that (\ref{StrMoRi1}) remains valid for $i=1,2,3.$
It is simple to check that $\displaystyle\sum_{i=1}^{3}2/q_i=2$,  $\displaystyle\sum_{i=1}^{3}
\delta(r_i)=3-2\displaystyle\sum_{i=1}^{3}1/r_i=1$, $\displaystyle\sum_{i=1}^{3}\theta_i=\frac{2}{1+\epsilon}$. Hence, adding the three equations in (\ref{StrBass6}), we see that necessarily, $\sum 2/q_i=\sum\delta(ri) +3-\sum\theta_i-\frac{\beta_1\delta(r_1)}{2}$. By the relation in (\ref{StrBass7}) This relation is equivalent  $2=4-\frac{2}{1+\epsilon}-\frac{b}{6}$ i.e.
\begin{equation}\label{StrBass8}
\frac{\epsilon}{1+\epsilon}=\frac{b}{12}
\end{equation}
Therefore , for $b$ small enough, it is clear that $\epsilon=\epsilon(b)=0^+$. Now, we choose
$(r_1,r_2,r_3)=(\frac{4}{1+b},\frac{2}{1-b},\frac{4}{1+b})$. It is simple to see that
$\sum1/r_i=1$ and $(\delta(r_1),\delta(r_2),\delta(r_3))=(1/2-b/2,b,1/2-b/2)$ and the relations of
(\ref{StrBass6})
can be written
\begin{equation}\label{StrBass9}
2/q_1=1/2-b/2+\frac{\epsilon +2b}{1+\epsilon},\hspace{0.2 cm} 2/q_2=b+1-\frac{2b}{1+\epsilon}-
\frac{b}{6},\hspace{0.2 cm}2/q_3=1/2-b/2+\frac{\epsilon}{1+\epsilon}.
\end{equation}
Now, by using the relation (\ref{StrBass8}), it results that $\epsilon=\frac{b}{12-b}$
and we see,
\begin{equation}\label{StrBass10}
2/q_1=1/2-b/2+b/12 +2b(1-b/12)=1/2+\frac{19b-2b^2}{12},
\end{equation}
\begin{equation}\label{StrBass11}
2/q_2=1- 5b/6-2b(1-b/12)=1+\frac{-7b+b^2}{6},
\end{equation}
\begin{equation}\label{StrBass12}
2/q_3=1/2-b/2+b/12=1/2-5b/12.
\end{equation}
therefore by the relations (\ref{StrBass10}), (\ref{StrBass11}), (\ref{StrBass12}), it is clear
 for $b$ small enough that $ (2/q_1,2/q_2,2/q_3)=(\frac{1}{2}^+,1^-, \frac{1}{2}^-)$ i.e.
 $(q_1,q_2,q_3)=(4^-,2^+,4^+)$  and by construction we have  $\sum1/q_i=1$. It remains
be checked that (\ref{StrMoRi1}) is valid for our parameter. Indeed, by definition of $\theta_i$ in
(\ref{StrBass5}) it is clear that $0<\theta_i<1$ for $i=1,2,3$. Since $(\delta(r_2),\theta_2)=
(b,\frac{2b}{1+\epsilon})$ we see that
\begin{equation}\label{StrBass13}
0\leq\beta_2=\frac{b}{2\delta(r_2)}=\frac{1}{2}, \hspace{0.2 cm} 0\leq \delta(r_2)=b\leq 2b(1-b/12)=
\frac{2b}{1+\epsilon}=\theta_2\leq\frac{\theta_2}{1-\beta/3}
\end{equation}
Moreover, since we have the equality $(\delta(r_1),\delta(r_2))=(\frac{1}{2}^-,\frac{1}{2}^-)$ and
$(\theta_1,\theta_2)=(\frac{1-2b}{1+\epsilon},\frac{1}{1+\epsilon})$  ,it is simple to see for $i=1,2$  that,
\begin{equation}\label{StrBass14}
0\leq\beta_i\leq 1/2, \hspace{0.2 cm} 0\leq \delta(r_i)\leq \theta_i\leq \frac{\theta_i}{1-\beta_i/3}
\end{equation}
By combining the equations  (\ref{StrBass13}) and (\ref{StrBass14}), we have (\ref{StrMoRi1}) and
 now we can apply  Lemma \ref{StrMoRi} to obtain
\begin{equation}
||F^{-1}_{t,x}(\frac{|\hat w(\tau,\nu)|}{\langle \sigma \rangle^{1/2-b}})||_{L^{q_1,r_1}_{t,x}}\leq
CT^\mu||w||_{L^2_{\tau,\nu}}
\end{equation}
\begin{equation}
||F^{-1}_{t,x}(\frac{|\hat u(\tau,\nu)|}{\langle \sigma \rangle^{b}|\xi|^{b/4}})||_
{L^{q_2,r_2}_{t,x}}\leq CT^\mu||u||_{L^2_{\tau,\nu}}
\end{equation}
\begin{equation}
||F^{-1}_{t,x}(\frac{|\hat v(\tau,\nu)|}{\langle \sigma \rangle^{1/2}})||_{L^{q_3,r_3}_{t,x}}
\leq C||v||_{L^2_{\tau,\nu}}
\end{equation}
 This completes the proof.$\hfill{\Box}$\\
\begin{lemma}\label{Str1MoRi}
Let $u$, $v\in L^2(\mathbb R^3)$ with compact support in $\{(x,y,t):|t|\leq T\}$. For $b>0$ and  $c>0$   small enough  there exists $\mu >0$ such that

\begin{eqnarray}\label{Str1MoRi2}
\nonumber && \int_{\mathbb R^6}\frac{|\hat u(\tau_1,\nu_1)||\hat v(\tau-\tau_1,\nu-\nu_1)||\hat w(\tau,\nu)|}
{\langle \sigma_1 \rangle^{1/2}|\xi_1|^{3b+c}\langle\sigma_2\rangle^{1/2-b}}d\tau d\tau_1 d\nu d\nu_1\\
&&\hspace{2cm}\leq CT^\mu ||u||_{L^2_{t,x}}||v||_{L^2_{t,x}}||w||_{L^2_{t,x}}
\end{eqnarray}
where $\sigma$, $\sigma_1$ and $\sigma_2$ are defined by
\begin{equation}
\sigma=\tau-P(\nu),\hspace{0,2cm} \sigma_1=\tau_1-P(\nu_1), \hspace{0,2cm} \sigma_2=\tau-\tau_1 -
P(\nu-\nu_1)
\end{equation}$\hfill{\Box}$
\end{lemma}
{\bf Proof.} By Plancherel theorem and by H\"older inequatliy  in space and time we see that the
 right hand side of (\ref{Str1MoRi2}) is bounded by
\begin{equation}\label{SMRB1}
||\mathcal F^{-1}_{t,x}(\frac{|\hat u(\tau,\nu)|}{\langle \sigma \rangle^{1/2}|\xi|^{3b+c}})||_{L^{q_1,r_1}_{t,\nu}}
||\mathcal F^{-1}_{t,x}(\frac{|\hat v(\tau,\nu)|}{\langle \sigma \rangle^{1/2-b}})||_{L^{q_2,r_2}_{t,\nu}}
||w||_{L^2_{t,\nu}}
\end{equation}
provided
\begin{equation}\label{SMRB2}
1/r_1 +1/r_2=1/2,\quad 1/q_1+1/q_2=1/2
\end{equation}
To apply Lemma \ref{StrMoRi} to each of the first two terms in (\ref{SMRB1}), for $b$ and $c$
small enough we take $\epsilon_1=\epsilon_2=\epsilon $, where $\epsilon=\epsilon(b,c)$
 We set,
\begin{equation}\label{SMRB3}
\theta_1=\frac{1}{1+\epsilon},\quad \theta_2=\frac{1-2b}{1+\epsilon}
\end{equation}
we choose $\beta_2=0$ and $\beta_1$ such that  $\frac{\beta_1\delta(r_1)}{2}= 3b - c$.   From
(\ref{StrMoRi3}), it remains to find $\beta_1$, $q_i$ and $r_i$ with
\begin{equation}\label{SMRB4}
\frac{2}{q_1}=\delta(r_1)+(1-\theta_1)-\frac{\beta_1\delta(r_1)}{3},\quad \frac{2}{q_2}=\delta(r_2)+(1-\theta_2)
\end{equation}
such that (\ref{StrMoRi1}) remains valid for $i=1,2.$
It is simple to see that $\displaystyle\sum_{i=1}^{2}2/q_i=1$,  $\displaystyle\sum_{i=1}^{2}
\delta(r_i)=1$ and  $\displaystyle\sum_{i=1}^{3}\theta_i=\frac{2-2b}{1+\epsilon}$. Hence, adding the two  equations of (\ref{SMRB4}), we see that necessarily,$ 1= 1 + (2-\frac{2-2b}{1+\epsilon})-
\frac{6b+2c}{3}$,   This relation is equivalent  $\frac{3-3b-2c}{3}=\frac{1-b}{1+\epsilon}$ i.e.
\begin{equation}\label{SMRB5}
\epsilon=\frac{c}{3 - 3b - c}
\end{equation}
Therefore,  for  $b$  and $ c$ small enough, it is clear that $\epsilon=\epsilon(b,c)=0^+$. Now  we choose
$(r_1,r_2)=(4,4)$ It follows that
$\sum1/r_i=1$ and $(\delta(r_1),\delta(r_2))=(1/2,1/2)$ and the relations of
(\ref{StrBass6}) can be written
\begin{equation}\label{SMRB6}
2/q_1=1/2 +\frac{\epsilon}{1+\epsilon} - \frac{6b+2c}{3},\hspace{0.2 cm} 2/q_2= 1/2 +\frac{2b+\epsilon}{1+\epsilon}.
\end{equation}
 From  (\ref{SMRB5}),  we get that $(2/q_1,2/q_2)=({\frac{1}{2}}^-,{\frac{1}{2}}^+)$ i.e.
 $(q_1,q_2,)=(4^+,4^-)$  and by construction we have  $\sum1/q_i=1/2$.
Moreover  (\ref{StrMoRi1}) is valid for our parameter. Indeed, for $b$ and $c$ small we have that
$$(\theta_1 ,\theta_2)=(1^-,1^-),\quad \beta_2=\frac{6b+2c}{\delta(r_2)}= 12b+4c=0^+\leq1/2$$
and
$$\delta(r_i)\sim 1/2 <1^-=\theta_i\leq\frac{\theta_i}{1-\beta_i/3}$$
Now we  apply Lemma \ref{StrMoRi} to give un suitable  bound for each of the first  two terms in
(\ref{SMRB1})
 This ends  the proof of Lemma \ref{Str1MoRi}. $\hfill{\Box}$\\

\begin{lemma}\label{StrTzev}
Let $2\leq q\leq 4$ and $u\in L^2(\mathbb R^2)$ with compact support in $\{(x,y,t):|t|\leq T\}$. For
 $\epsilon>0$ and  $b=(1-2/q)(\frac{1+\epsilon}{2}),$  there exists $\mu=\mu(\epsilon)>0$  such that
\begin{equation}\label{StrTzev1}
||\mathcal F^{-1}(\langle \sigma \rangle^{-b}|\hat u(\tau, \nu)|)||_{L^q_{t,\nu}}\leq CT^\mu
 ||u||_{L^2_{t,\nu}}
\end{equation}$\hfill{\Box}$
\end{lemma}
$\bold{Proof}$. For any $\phi \in L^2(\mathbb R^2)$ the Strichartz inequality in \cite{Sau93}
(see Proposition 2.3) yields
\begin{equation}\label{StrTzev2}
||U(t)\phi||_{L^4_{t,\nu}}\leq ||\phi||_{L^2_\nu}
\end{equation}
 where  $$\mathcal F_{t,\nu}( U(t)\phi)=\exp(it(\xi^3 + n\frac{\eta^2}{\xi}))\hat \phi(\xi,\eta),
\quad n=\pm 1,$$
 using (\ref{StrTzev2}) together with Lemma 3.3 of \cite{Ge}, we see for all $\epsilon >0$ that
\begin{equation}\label{StrTzev3}
|| u(\tau, \nu)||_{L^4_{t,\nu}}\leq C
 ||\langle \sigma \rangle^{1/2 +\epsilon/4}\hat u(\tau,\nu)||_{L^2_{t,\nu}}.
\end{equation}

Since we have $ ||u||_{L^2_{t,\nu}}=||u||_{X^{0,0,0}}$, therefore by interpolation between
$(L^4_{t,\nu}, X^{1/2+\epsilon/4,0,0})$ and $ (L^2_{t,\nu}, X^{0,0,0})$, for $0\leq \theta \leq 1$, we obtain
 \begin{equation}\label{StrTzev31}
|| u(\tau, \nu)||_{L^q_{t,\nu}}\leq C
 ||\langle \sigma \rangle^{\theta(1/2 +\epsilon/4)}\hat u(\tau,\nu)||_{L^2_{t,\nu}}
\end{equation}
Where
\begin{equation}\label{StrTzev4}
2/q=\theta/4 + \frac{1-\theta}{2}
\end{equation}
Next, using the assumption on the support of u and the results in \cite{Ge95}, we get that there
exists $\mu=\mu(\epsilon)$ such that
\begin{equation}\label{StrTzev5}
|| u(\tau, \nu)||_{L^q_{t,\nu}}\leq C T^\mu
 ||\langle \sigma \rangle^{\theta(1/2 +\epsilon/2)}\hat u(\tau,\nu)||_{L^2_{t,\nu}}
\end{equation}
 from (\ref{StrTzev4}), the desired result is deduced. $\hfill{\Box}$\\
Using Lemma \ref{StrTzev}  together whith the Proof of Lemma 2.2 of  \cite{TaTz}, we obtain the
following Lemma.
\begin{lemma}\label{StrTzBas}
Let $u$, $v$,  $w\in L^2(\mathbb R^3)$ with compact support in $\{(x,y,t):|t|\leq T\}$ and  $\alpha$,
$\beta$, $\gamma\in [0,1/2+\epsilon]$ . For any $\epsilon>0$ there exists $\mu=\mu(\epsilon) >0$ such that
\begin{eqnarray}\label{StrTzBas1}
\nonumber&& \int_{\mathbb R^6}\frac{|\hat u(\tau_1,\nu_1)||\hat v(\tau-\tau_1,\nu-\nu_1)||\hat w(\tau,\nu)|}
{\langle \sigma \rangle^{\alpha}\langle\sigma_1\rangle^{\beta} \langle\sigma_2\rangle^{\gamma}}d\tau d\tau_1 d\nu d\nu_1\\
&&\hspace{2cm}\leq CT^\mu ||u||_{L^2_{t,x}}||v||_{L^2_{t,x}}||w||_{L^2_{t,x}},
\end{eqnarray}
proven for $\alpha+\beta+\gamma \geq 1+2\epsilon$ .$\hfill{\Box}$
\end{lemma}
\section{Bilinear estimates}\label{EstBL}
In this section we will prepare certain bilinear estimates on  $\partial_x(uv)$ in the Bourgain
space $X^{b,s_1,s_2}$. These bilinear estimates will be the main tools in the next section to apply a fixed point argument which will give the local existence result.
\begin{proposition}\label{EstBlBas}
Let  $\delta>0$  small enough, $s_2\geq 0$ and $s_1\in [\frac{-1}{2}+8\delta,0]$. For all $u$, $v\in X^{1/2,s_1,s_2}$  whith compact support in time and  included in the subset $\{(t,x,y):t\in[-T,T]\}$,
there exists $\mu>0$ such that the following  bilinear estimate holds
\begin{equation}\label{EstBlBas1}
||\partial_x(uv)||_{X^{-1/2+\delta,s_1-2\delta +\epsilon,s_2}}\leq CT^\mu ||u||_{X^{1/2,s_1,s_2}}
||v||_{X^{1/2,s_1,s_2}}
\end{equation}
for some  $\epsilon>0$ such that $\epsilon<<\delta$.$\hfill{\Box}$
\end{proposition}
$\bold{Proof}$. We proceed  by duality.  It is equivalent to show that for
 $\delta >0$ small enough and $\epsilon<<\delta$ for all $w\in X^{1/2-\delta,-s_1+2\delta-\epsilon,-s_2}$,
 \begin{equation}\label{EstBlBas2}
|\langle\partial_x(uv),w\rangle|\leq CT^\mu [||u||_{X^{1/2,s_1,s_2}}||v||_{X^{1/2,s_1,s_2}}]
||w||_{ X^{1/2-\delta,-s_1+2\delta -\epsilon,-s_2}}
\end{equation}
Let $f$, $g$ and $h$ respectively defined by
  \begin{equation}\label{EstBlBas3}
\hat f(\tau,\nu)=\langle i(\tau-P(\nu)) +\xi^2\rangle^{1/2}\langle \xi\rangle^{s_1}
\langle \eta \rangle ^{s_2} \hat u(\tau, \nu),
\end{equation}
 \begin{equation}\label{EstBlBas4}
\hat g(\tau,\nu)=\langle i(\tau-P(\nu)) +\xi^2\rangle^{1/2}\langle \xi\rangle^{s_1}
\langle \eta \rangle ^{s_2} \hat v(\tau, \nu),
\end{equation}
 \begin{equation}\label{EstBlBas5}
\hat h(\tau,\nu)=\langle i(\tau-P(\nu)) +\xi^2\rangle^{-1/2+\delta}\langle \xi\rangle^{-s_1+2\delta-\epsilon}
\langle \eta \rangle ^{-s_2} \hat w(\tau, \nu).
\end{equation}
It is clear that
$$||u||_{X^{1/2,s_1,s_2}}=||f||_{L^2_{t,\nu}},\hspace{0.1cm}
||v||_{X^{1/2,s_1,s_2}}=||g||_{L^2_{t,\nu}},\hspace{0.1cm}
||w||_{X^{-1/2+\delta,-s_1+2\delta -\epsilon,-s_2}}=||h||_{L^2_{t,\nu}}$$

Thus  by Plancherel Theorem,(\ref{EstBlBas2})  is equivalent
 \begin{eqnarray}\label{EstBlBas6}
\nonumber&& \int_{\mathbb R^6}\frac{|\xi||\hat h(\tau_1,\nu_1)||\hat g(\tau-\tau_1,\nu-\nu_1)
||\hat h(\tau,\nu)|}{\langle i\sigma+\xi^2 \rangle^{1/2-\delta}
\langle i\sigma_1+\xi_1^2\rangle^{1/2} \langle i\sigma_2 +(\xi-\xi_1)^2\rangle^{1/2}}
  \frac{\langle \xi \rangle ^{s_1 -2\delta+\epsilon}}{\langle \xi -\xi_1\rangle^{s_1}
\langle \xi_1\rangle^{s_1}}\\
\nonumber&&\hspace{3cm}\times   \frac{\langle \eta \rangle ^{s_2}}{\langle \eta -\eta_1\rangle^{s_2}
\langle \eta_1\rangle^{s_2}}      d\tau d\tau_1 d\nu d\nu_1\\
&&\hspace{2cm}\leq CT^\mu ||u||_{L^2_{t,x}}||v||_{L^2_{t,x}}||w||_{L^2_{t,x}}
\end{eqnarray}

 Moreover, we can assume that $s_2=0$ since in the case  $s_2\geq 0$
we have
\begin{equation}\label{EstBlBas7}
 \frac{\langle \eta \rangle ^{s_2}}{\langle \eta -\eta_1\rangle^{s_2}
\langle \eta_1\rangle^{s_2}}\leq C, \quad \forall \eta,\hspace{0.1 cm} \eta_1\in \mathbb R.
\end{equation}
Therefore,setting $s=- s_1\in[0,1/2-8\delta]$, it is enough to estimate,
\begin{eqnarray}\label{EstBlBas8}
\nonumber I&=:& \int_{\mathbb R^6}\frac{|\hat h(\tau_1,\nu_1)||\hat g(\tau-\tau_1,\nu-\nu_1)
||\hat h(\tau,\nu)|}{\langle i\sigma+\xi^2 \rangle^{1/2-\delta}
\langle i\sigma_1+\xi_1^2\rangle^{1/2} \langle i\sigma_2 +(\xi-\xi_1)^2\rangle^{1/2}}\\
 &&\hspace{2cm}\times \frac{|\xi|\langle \xi -\xi_1\rangle^s\langle \xi_1\rangle^s}
{\langle \xi\rangle^{s+2\delta-\epsilon}} d\tau d\tau_1 d\nu d\nu_1
\end{eqnarray}
To estimate $I$ we will use a algebric relation  between $\sigma$, $\sigma_1$ and $\sigma_2$:
\begin{equation}\label{EstBlBas9}
\sigma_1+\sigma_2-\sigma= 3\xi \xi_1(\xi-\xi_1) +\frac{(\xi_1\eta -\xi \eta_1)^2}{\xi \xi_1(\xi-\xi_1)}
\end{equation}
(see \cite{Bo93}), which ensures that
\begin{equation}\label{EstBlBas10}
\max(|\sigma|,|\sigma_1|,|\sigma_2|)\geq \frac{|\sigma_1+\sigma_2-\sigma|}{3}\geq
 |\xi \xi_1(\xi-\xi_1)|
\end{equation}
A symmetry argument  shows that it is enough to estimate the contribution to  $I$ of the
subset  of $\mathbb R^6$, $\Omega=\{(\tau,\tau_1,\nu,\nu_1)\in\mathbb R^6:|\sigma_1|\geq|\sigma_2|\}
$. To do this  we split $\Omega$ in $\Omega=\Omega_1\cup \Omega_2$ where
$$ \Omega_1=\Omega\cap\{(\tau,\tau_1,\nu,\nu_1)\in\mathbb R^6: |\xi|\leq C_0, C_0>>1\},$$
$$ \Omega_2=\Omega\cap\{(\tau,\tau_1,\nu,\nu_1)\in\mathbb R^6: |\xi|\geq C_0, C_0>>1\},$$
{\bf Case 1}~: Contribution of  $\Omega_1$ to $I$. We  divide $\Omega_1$
in three region:
 $$ \Omega_1^1=\Omega_1\cap\{(\tau,\tau_1,\nu,\nu_1)\in\mathbb R^6: |\xi_1|\leq 2C_0\},$$
 $$ \Omega_1^2=\Omega_1\cap\{(\tau,\tau_1,\nu,\nu_1)\in\mathbb R^6:|\sigma|\geq|\sigma_1|,
 |\xi_1|\geq 2C_0\},$$
 $$ \Omega_1^3=\Omega_1\cap\{(\tau,\tau_1,\nu,\nu_1)\in\mathbb R^6:|\sigma_1|\geq|\sigma|,
 |\xi_1|\geq 2C_0\}.$$
It is clear that $\Omega_1=\Omega_1^1\cup \Omega_2^2\cup\Omega_1^3$.\\
{ \bf Case 1.1}~: Contribution of  $\Omega_1^1$ to $I$. Denote by  $I_1^1$ the contribution of this region to $I$. In this case we have $|\xi -\xi_1|\leq |\xi|+|\xi_1|\leq C$ and we see that
$$ \frac{|\xi|\langle \xi -\xi_1\rangle^s\langle \xi_1\rangle^s}
{\langle \xi\rangle^{s+2\delta-\epsilon}}\leq C$$
 and hence,
\begin{eqnarray*}
 I_1^1 &\leq&C \int_{\mathbb R^6}\frac{|\hat h(\tau_1,\nu_1)||\hat g(\tau-\tau_1,\nu-\nu_1)
||\hat h(\tau,\nu)|}{\langle i\sigma+\xi^2 \rangle^{1/2-\delta}
\langle i\sigma_1+\xi_1^2\rangle^{1/2} \langle i\sigma_2 +(\xi-\xi_1)^2\rangle^{1/2}}
 d\tau d\tau_1 d\nu d\nu_1\\
&\leq& C \int_{\mathbb R^6}\frac{|\hat h(\tau_1,\nu_1)||\hat g(\tau-\tau_1,\nu-\nu_1)
||\hat h(\tau,\nu)|}{\langle \sigma \rangle^{1/2-\delta}
\langle \sigma_1\rangle^{1/2} \langle \sigma_2\rangle^{1/2}}d\tau d\tau_1 d\nu d\nu_1
\end{eqnarray*}
Now, we can apply  Lemma \ref{StrTzBas} to deduce
$$I_1^1\leq CT^\mu ||u||_{L^2_{t,x}}||v||_{L^2_{t,x}}||w||_{L^2_{t,x}}$$
{ \bf Case 1.2}~: Contribution  of $\Omega_1^2$ to $I$. Denote by  $I_1^1$ the contribution of
this region to $I$. Since we have, In this case, $|\xi|\leq 1/2 |\xi_1|$, it  follows  that
$ |\xi_1|\sim |\xi-\xi_1|$. Therefore
   $$ \frac{|\xi|\langle \xi -\xi_1\rangle^s\langle \xi_1\rangle^s}
{\langle \xi\rangle^{s+2\delta -\epsilon}}\leq C |\xi-\xi_1|^{2s}|\xi|^s.$$
  Moreover, since $|\sigma|=\max(|\sigma|,|\sigma_1|,|\sigma_2|)$ , by the
  relation  between $\sigma$, $\sigma_1$ and $\sigma_2$ in (\ref{EstBlBas9})  it
 results that $|\sigma|\geq |\xi||\xi_1||\xi-\xi_1|$. Therefore
$$|\sigma|^s\geq |\xi|^s|\xi_1||\xi-\xi_1|\geq |\xi-\xi_1|^{2s}|\xi|^s$$
and hence,
 \begin{eqnarray*}
 I_1^2 &\leq&C \int_{\mathbb R^6}\frac{|\xi-\xi_1|^{2s}|\xi|^s|\hat h(\tau_1,\nu_1)||\hat
 g(\tau-\tau_1,\nu-\nu_1)
||\hat h(\tau,\nu)|}{\langle i\sigma+\xi^2 \rangle^{1/2-\delta}
\langle i\sigma_1+\xi_1^2\rangle^{1/2} \langle i\sigma_2 +(\xi-\xi_1)^2\rangle^{1/2}}
 d\tau d\tau_1 d\nu d\nu_1\\
&\leq& C \int_{\mathbb R^6}\frac{|\hat h(\tau_1,\nu_1)||\hat g(\tau-\tau_1,\nu-\nu_1)
||\hat h(\tau,\nu)|}{\langle \sigma \rangle^{1/2-s-\delta}
\langle \sigma_1\rangle^{1/2} \langle \sigma_2\rangle^{1/2}}d\tau d\tau_1 d\nu d\nu_1,
\end{eqnarray*}
  since $s\in [0,1/2-8\delta]$, we see that $(1/2-s-\delta)+1/2+1/2> 1+8\delta$, therefore a use
 of Lemma \ref{StrTzBas}  provides a good bound  for $I_1^1$.\\\\
  {\bf Case 1.3}~: Contribution  of  $\Omega_1^3$ to $I$.  We denote  by  $I_1^3$ the
 contribution
 of this region to $I$. In this case $|\sigma_1|$ dominates and $ |\xi_1|\sim |\xi-\xi_1|$.
 Because of (\ref{EstBlBas9}) , we obtain
$$ \frac{|\xi|\langle \xi -\xi_1\rangle^s\langle \xi_1\rangle^s}{\langle \xi\rangle^{s+2\delta-\epsilon}}
\leq C |\xi-\xi_1|^{2s}|\xi|^s\leq \langle\sigma_1 \rangle^s.$$
As in case 1.2, we obtain by the Lemma  \ref{StrTzBas} that
 \begin{eqnarray*}
I_1^3 &\leq&C\int_{\mathbb R^6}\frac{|\hat h(\tau_1,\nu_1)||\hat g(\tau-\tau_1,\nu-\nu_1)
||\hat h(\tau,\nu)|}{\langle \sigma \rangle^{1/2-\delta}
\langle \sigma_1\rangle^{8\delta} \langle \sigma_2\rangle^{1/2}}d\tau d\tau_1 d\nu d\nu_1\\
&\leq& CT^\mu ||f||_{L^2_{t,x}}||g||_{L^2_{t,x}}||h||_{L^2_{t,x}}.
\end{eqnarray*}
 { \bf Case 2}~: Contribution  of  $\Omega_2$ to $I$. We divide $\Omega_2$ into three subdomain
$\Omega_2^i$,  $ i=1,2,3$ such that $\Omega_2=\Omega_2^1\cup\Omega_2^2\cup  \Omega_2^3$, where
 $$ \Omega_2^1=\Omega_2\cap\{(\tau,\tau_1,\nu,\nu_1)\in\mathbb R^6: \min(|\xi_1|,|\xi_2|)\leq 1\},$$
 $$ \Omega_2^2=\Omega_2\cap\{(\tau,\tau_1,\nu,\nu_1)\in\mathbb R^6:|\sigma|\geq|\sigma_1|,
 \min(|\xi_1|,|\xi_2|)\geq 1 \},$$
 $$ \Omega_2^3=\Omega_2\cap\{(\tau,\tau_1,\nu,\nu_1)\in\mathbb R^6:|\sigma_1|\geq|\sigma|,
 \min(|\xi_1|,|\xi_2|)\geq 1 \}.$$
{\bf Case 2.1}~: Contribution  de $I$ in $\Omega_2^1$.  We denote   by  $I_2^1$ the contribution
 of this region to $I$. By symmetry we can assume that $\min(|\xi_1|,|\xi-\xi_1|)=|\xi_1|$ and thus
$|\xi-\xi_1|\leq 1+ |\xi|\leq (1+C_0)|\xi|$ , therefore $|\xi|\sim|\xi-\xi_1|$ . It follows
$$ \frac{|\xi|\langle \xi -\xi_1\rangle^s\langle \xi_1\rangle^s}{\langle \xi\rangle^{s+2\delta-\epsilon}}\leq C |\xi|^{1-2\delta+\epsilon}.$$
Since $\langle i\sigma +\xi^2\rangle^{1/2-\delta}\geq |\xi|^{1-2\delta}$,  $\langle i\sigma_2+(\xi-\xi_2)^2\rangle^{1/2}\geq \langle\sigma_2\rangle^{1/2-\delta}|\xi-\xi_1|^{\epsilon}$ and $|\xi_1|\leq 1$,
 it results that
 \begin{eqnarray*}
I_2^1 &\leq&C\int_{\mathbb R^6}\frac{|\hat h(\tau_1,\nu_1)||\hat g(\tau-\tau_1,\nu-\nu_1)
||\hat h(\tau,\nu)|}{\langle \sigma_1 \rangle^{1/2}
 \langle \sigma_2\rangle^{1/2-\delta}}d\tau d\tau_1 d\nu d\nu_1\\
&\leq&\int_{\mathbb R^6}\frac{|\hat h(\tau_1,\nu_1)||\hat g(\tau-\tau_1,\nu-\nu_1)
||\hat h(\tau,\nu)|}{\langle \sigma_1 \rangle^{1/2}|\xi_1|^{4\delta}
 \langle \sigma_2\rangle^{1/2-\delta}}d\tau d\tau_1 d\nu d\nu_1.
\end{eqnarray*}
Now a use  of  Lemma \ref{Str1MoRi} provides a bound for $I_2^1$.

{\bf Case 2.2}~: Contribution  of  $\Omega_2^2$ to $I$. In this case  we need to divide  $\Omega_2^2$  in two  regions  defined by,
 $$\Omega_2^{21}=\Omega_2^2\cap\{(\tau,\tau_1,\nu,\nu_1)\in\mathbb R^6:
 \min(|\xi_1|,|\xi_2|)\leq \frac{1}{C_0} |\xi|\},$$
 $$\Omega_2^{22}=\Omega_2^2\cap\{(\tau,\tau_1,\nu,\nu_1)\in\mathbb R^6:
 \min(|\xi_1|,|\xi_2|)\geq \ \frac{1}{C_0} |\xi|\},$$
{\bf Case 2.21}: Contribution  of $\Omega_2^{21}$ to $I$.  We denote  by  $I_2^{21}$
the contribution of this region to $I$. By symmetry argument we can assume that
$|\xi_1|\leq |\xi-\xi_1| $, It follows $|\xi|\leq\frac{1}{C_0} |\xi|$, ($C_0>> 1$). Therefore
$|\xi -\xi_1|\leq |\xi|+|\xi-1|\leq C|\xi|$ and $|\xi|\leq |\xi_1| +|\xi -\xi_1|\leq
 \frac{1}{C_0} |\xi| +|\xi -\xi_1| $ i.e. $|\xi|\leq \frac{1}{1-1/C_0}|\xi-\xi_1|$ and thus
$|\xi-\xi_1|\sim|\xi|$. It results that
$$ \frac{|\xi|\langle \xi -\xi_1\rangle^s\langle \xi_1\rangle^s}{\langle \xi\rangle^{s+2\delta-\epsilon}}
\leq C |\xi|^{1-2\delta+\epsilon}|\xi_1|^s,$$
 and hence
$$I_2^{21} \leq C \int_{\mathbb R^6}\frac{|\xi|^{1-2\delta+\epsilon}|\xi_1|^s|\hat h(\tau_1,\nu_1)||\hat g(\tau-\tau_1,\nu-\nu_1)||\hat h(\tau,\nu)|}{\langle i\sigma+\xi^2 \rangle^{1/2-\delta}
\langle i\sigma_1+\xi_1^2\rangle^{1/2} \langle i\sigma_2 +(\xi-\xi_1)^2\rangle^{1/2-\delta}|\xi-\xi_1|^\epsilon}|
 d\tau d\tau_1 d\nu d\nu_1.$$
Since in this case $|\sigma|$ dominates, we obtain
$$\langle i \sigma +\xi^2\rangle^{1/2-\delta}\geq \langle \sigma \rangle^{1/2-\delta}\geq
|\xi|^{1/2-\delta}|\xi_1|^{1/2-\delta}|\xi-\xi_1|^{1/2-\delta}\geq |\xi|^{1-2\delta}
|\xi_1|^{1/2-\delta}$$
therefore,
$$I_2^{21} \leq C \int_{\mathbb R^6}\frac{|\hat h(\tau_1,\nu_1)||\hat g(\tau-\tau_1,\nu-\nu_1)||\hat h(\tau,\nu)|}{|\xi_1|^{1/2-\delta-s}
\langle \sigma_1\rangle^{1/2} \langle \sigma_2 \rangle^{1/2-\delta}}
 d\tau d\tau_1 d\nu d\nu_1.$$
It is clear that $1/2-\delta-s\geq 4\delta$ for $s\in[0,1/2-8\delta]$. Now we can apply
  Lemma \ref{Str1MoRi} to estimate  $I_2^{21}$.\\\\
{\bf Case 2.22}~: Contribution  of $\Omega_2^{22}$ to $I$.  We denote  by  $I_2^{22}$
the contribution of this region to $I$. In this case, we notice that $|\xi|\lesssim |\xi_1|$
and    $|\xi|\lesssim |\xi-\xi_1|$, it results that:
 \begin{eqnarray*}
 \frac{|\xi|\langle \xi -\xi_1\rangle^s\langle \xi_1\rangle^s}{\langle \xi\rangle^{s+2\delta-\epsilon}}
&\leq& C |\xi|^{1-2\delta-s+\epsilon}|\xi_1|^s|\xi-\xi_1|^s\\
&\leq& C |\xi|^{\frac{1-2\delta-s+\epsilon}{3}}|\xi|^{\frac23(1-2\delta +\epsilon-s)}|\xi_1|^s|\xi-\xi_1|^s\\
&\leq& C |\xi|^{\frac{1-2\delta +\epsilon-s}{3}}|\xi_1|^{\frac{1-2\delta +\epsilon-s}{3}}
|\xi-\xi_1|^{\frac{1-2\delta +\epsilon -s}{3}}\\
&\leq&\langle \sigma \rangle ^{\frac{1-2\delta +\epsilon -s}{3}},
\end{eqnarray*}
and hence
$$I_2^{22} \leq C \int_{\mathbb R^6}\frac{|\hat h(\tau_1,\nu_1)||\hat g(\tau-\tau_1,\nu-\nu_1)||\hat h(\tau,\nu)|}{\langle \sigma \rangle ^{1/2-\frac{1-2\delta +\epsilon-s}{3}-\delta}
\langle \sigma_1\rangle^{1/2} \langle \sigma_2 \rangle^{1/2}}
 d\tau d\tau_1 d\nu d\nu_1.$$
since $s\in [0,1/2-8\delta]$ we see, for $\epsilon<<\delta$, that
 $ 1/2-\frac{1-2\delta+\epsilon-s}{3}-\delta=\frac{1-2s+4\delta-2\epsilon}{3}\geq 2\delta$, therefore we can apply
 Lemma \ref{StrTzBas} to estimate $I_2^{22}$.\\\\
{\bf Case 2.3}~: Contribution  of  $\Omega_2^3$ to $I$. We  divide $\Omega_2^3$ in two parts:
  $$\Omega_2^{31}=\Omega_2^3\cap\{(\tau,\tau_1,\nu,\nu_1)\in\mathbb R^6:
 \min(|\xi_1|,|\xi_2|)\leq \frac{1}{C_0} |\xi|\},$$
 $$\Omega_2^{32}=\Omega_2^3\cap\{(\tau,\tau_1,\nu,\nu_1)\in\mathbb R^6:
 \min(|\xi_1|,|\xi_2|)\geq \ \frac{1}{C_0} |\xi|\},$$
 {\bf Case 2.31}~: Contribution  of  $\Omega_2^{31}$ to $I$. Because there is no symmetry between  $|\xi_1|$ and $|\xi-\xi_1|$ we distinguish between two regions of  $\Omega_2^{31}$~:
   $$\Omega_2^{311}=\Omega_2^{31}\cap\{(\tau,\tau_1,\nu,\nu_1)\in\mathbb R^6:
 \min(|\xi_1|,|\xi_2|)=|\xi_1|\},$$
 $$\Omega_2^{312}=\Omega_2^{31}\cap\{(\tau,\tau_1,\nu,\nu_1)\in\mathbb R^6:
 \min(|\xi_1|,|\xi_2|)=|\xi-\xi_1|\},$$
{\bf Case 2.311}~: Contribution of $\Omega_2^{311}$ to $I$.  We denote  by  $I_2^{311}$
the contribution of this region to $I$. In this case  we have $|\xi_1|\leq \frac{1}{C_0} |\xi|$
 and thus $|\xi-\xi_1|\sim |\xi|$. Therefore
\begin{eqnarray*}
 \frac{|\xi|\langle \xi -\xi_1\rangle^s\langle \xi_1\rangle^s}{\langle \xi\rangle^{s+2\delta-\epsilon}}
&\leq C&|\xi|^{1-2\delta+\epsilon}||\xi_1|^s\\
&\leq C&|\xi-\xi_1|^{1-2\delta+\epsilon}||\xi_1|^s
\end{eqnarray*}
it results that
$$I_2^{311} \leq C \int_{\mathbb R^6}\frac{|\xi-\xi_1|^{1-2\delta+\epsilon}|\xi_1|^s|\hat h(\tau_1,\nu_1)||\hat g(\tau-\tau_1,\nu-\nu_1)||\hat h(\tau,\nu)|}{\langle i\sigma+\xi^2 \rangle^{1/2-\delta}
\langle i\sigma_1+\xi_1^2\rangle^{1/2} \langle i\sigma_2 +(\xi-\xi_1)^2\rangle^{1/2}}
 d\tau d\tau_1 d\nu d\nu_1.$$
Since $|\sigma_1|$ dominates , we see that
\begin{eqnarray*}
\langle i \sigma_1 +\xi_1^2\rangle^{1/2}&\geq& \langle \sigma_1
\rangle^{\delta-\epsilon/2}\langle \sigma
\rangle^{\epsilon/2}\langle \sigma_1 \rangle^{1/2-\delta+\epsilon}\\
&\geq&\langle \sigma \rangle^{\epsilon/2} \langle \sigma_1 \rangle^{\delta-\epsilon/2}|\xi|^{1/2-\delta+\epsilon/2}|\xi_1|^{1/2-\delta+\epsilon/2}
|\xi-\xi_1|^{1/2-\delta+\epsilon/2}\\
&\geq& \langle \sigma \rangle^{\epsilon/2} \langle \sigma_1 \rangle^{\delta-\epsilon}|\xi-\xi_1|^{1-2\delta+\epsilon}|\xi_1|^{1/2-\delta}
\end{eqnarray*}
therefore
$$I_2^{311} \leq C \int_{\mathbb R^6}\frac{|\hat h(\tau_1,\nu_1)||\hat g(\tau-\tau_1,\nu-\nu_1)||\hat h(\tau,\nu)|}{\langle \sigma \rangle^{1/2-\delta+\epsilon/2}
\langle \sigma_1\rangle^{\delta-\epsilon/2}|\xi_1|^{1/2-\delta-s} \langle \sigma_2 \rangle^{1/2}}
 d\tau d\tau_1 d\nu d\nu_1.$$
Since  $1/2-\delta-s>\frac{\delta-\epsilon/2}{4}$, it follows  by virtue of Lemma \ref{StrBass}, for $\epsilon<<\delta$, that
\begin{eqnarray*}
I_2^{311} \leq &C& \int_{\mathbb R^6}\frac{|\hat h(\tau_1,\nu_1)||\hat g(\tau-\tau_1,\nu-\nu_1)||\hat h(\tau,\nu)|}{\langle \sigma \rangle^{1/2-\delta+\epsilon/2}
\langle \sigma_1\rangle^{\delta-\epsilon/2}|\xi_1|^{\frac{\delta-\epsilon/2}{4}} \langle \sigma_2 \rangle^{1/2}}
 d\tau d\tau_1 d\nu d\nu_1\\
&\leq& CT^\mu ||f||_{L^2_{t,x}}||g||_{L^2_{t,x}}||h||_{L^2_{t,x}}
\end{eqnarray*}
 {\bf Case 2.312}~: Contribution  of $\Omega_2^{312}$ to $I$.  We denote  by  $I_2^{312}$
the contribution of this region to $I$. In this case  we have $|\xi-\xi_1|\leq \frac{1}{C_0} |\xi|$
 and thus $|\xi_1|\sim |\xi|$. It results that
\begin{equation}·\label{EstBlBas11}
\frac{|\xi|\langle \xi -\xi_1\rangle^s\langle \xi_1\rangle^s}{\langle \xi\rangle^{s+2\delta-\epsilon}}
\leq C|\xi|^{1-2\delta+\epsilon}||\xi-\xi_1|^s.
\end{equation}
 Since $|\sigma_1|$ dominates in this case , for $\epsilon<\delta$ we obtain
\begin{eqnarray*}
\langle i \sigma_1 +\xi_1^2\rangle^{1/2}&\geq& \langle i\sigma_1+ \xi_1^2 \rangle^{\delta}\langle \sigma_1 \rangle^{1/2-\delta}\\
&\geq& |\xi|^\epsilon \langle \sigma_1 \rangle^{\delta-\epsilon/2}|\xi|^{1/2-\delta}|\xi_1|^{1/2-\delta}
|\xi-\xi_1|^{1/2-\delta}\\
&\geq&  \langle \sigma_1 \rangle^{\delta-\epsilon/2}|\xi-\xi_1|^{1/2-\delta}|\xi|^{1-2\delta+\epsilon}
\end{eqnarray*}
 and since $|\sigma_1|\geq|\sigma|$, it results, for $s\in [0,1/2-8\delta]$ and $\epsilon<\delta$,
 that
\begin{equation}\label{EstBlBas12}
\langle i \sigma_1 +\xi_1^2\rangle^{1/2}\geq
  \langle \sigma \rangle^{\delta-\epsilon/2}|\xi-\xi_1|^{1/2-\delta}|\xi|^{1-2\delta+\epsilon}.
\end{equation}
By combining (\ref{EstBlBas11}) and (\ref{EstBlBas12}), we can deduce that
\begin{eqnarray*}
I_2^{312} &\leq& C \int_{\mathbb R^6}\frac{|\hat h(\tau_1,\nu_1)||\hat g(\tau-\tau_1,\nu-\nu_1)|
|\hat h(\tau,\nu)|}{\langle \sigma \rangle^{1/2-\epsilon}
|\xi-\xi_1|^{1/2-s-\delta} \langle \sigma_2 \rangle^{1/2}}
 d\tau d\tau_1 d\nu d\nu_1\\
&\leq& C \int_{\mathbb R^6}\frac{|\hat h(\tau_1,\nu_1)||\hat g(\tau-\tau_1,\nu-\nu_1)|
|\hat h(\tau,\nu)|}{\langle \sigma \rangle^{1/2-\epsilon}
|\xi-\xi_1|^{4\epsilon} \langle \sigma_2 \rangle^{1/2}}
 d\tau d\tau_1 d\nu d\nu_1.
\end{eqnarray*}
  Now  Lemma \ref{Str1MoRi} provides a
 bound for $I_2^{312}$.\\\\
 {\bf Case 2.33}~: Contribution  of  $\Omega_2^{33}$ to $I$.  We indicate  by  $I_2^{33}$
the contribution of this region to $I$. In this case  $|\sigma_1|$ dominates and  we have
$|\xi-\xi_1|\geq \frac{1}{C_0} |\xi|$,  $|\xi_1|\geq \frac{1}{C_0} |\xi|$ . Hence
 \begin{eqnarray*}
 \frac{|\xi|\langle \xi -\xi_1\rangle^s\langle \xi_1\rangle^s}{\langle \xi\rangle^{s+2\delta-\epsilon}}
&\leq& C |\xi|^{1-2\delta+\epsilon-s}|\xi_1|^s|\xi-\xi_1|^s\\
&\leq& C |\xi|^{\frac{1-2\delta+\epsilon-s}{3}}|\xi_1|^{\frac{1-2\delta+\epsilon-s}{3}}
|\xi-\xi_1|^{\frac{1-2\delta+\epsilon-s}{3}}\\
&\leq&\langle \sigma_1 \rangle ^{\frac{1-2\delta+\epsilon-s}{3}},
\end{eqnarray*}
it follows
$$I_2^{33} \leq C \int_{\mathbb R^6}\frac{|\hat h(\tau_1,\nu_1)||\hat g(\tau-\tau_1,\nu-\nu_1)||\hat h(\tau,\nu)|}{\langle \sigma \rangle^{1/2-\delta}
\langle \sigma_1\rangle^{1/2-\frac{1-2\delta+\epsilon-s}{3}} \langle \sigma_2 \rangle^{1/2}}
 d\tau d\tau_1 d\nu d\nu_1.$$
it is clear that, for $\epsilon<<\delta$, $1/2-\frac{1-2\delta+\epsilon-s}{3}=\frac{1-2s+4\delta-2\epsilon}{3}\geq 3\delta$. therefore we can apply Lemma \ref{StrTzBas} to estimate $I_2^{33}$.\\
This completes the proof of theorem 1. $\hfill{\Box}$\\
Actually, we will mainly use the following version, which is a direct consequence of Proposition
 \ref{EstBlBas}, together with the  two  triangle inequality
\begin{equation}
\forall s_1\geq s_c^1, \quad \langle \xi\rangle^{s_1}\leq\langle\xi\rangle^{s_c^1}
\langle\xi_1\rangle^{s_1-s_c^1}+\langle\xi\rangle^{s_c^1}\langle\xi-\xi_1\rangle^{s_1-s_c^1},
\end{equation}
\begin{equation}
 \langle \eta \rangle^{s_2}\leq\langle\eta_1\rangle^{s_2}+
\langle\eta-\eta_1\rangle^{s_2}
\end{equation}
\begin{proposition}\label{EstBlBasM}
Let $s_c^1\in ]-1/2,0]$,  $s_2\geq 0$  . For all $s_1\geq s_c^1$ and  $ u$,
$v\in X^{1/2,s_1,s_2}$  whith compact support in time and  included in the subset
 $\{(t,x,y):t\in[-T,T]\}$,
there exists $\mu>0$ such that the following  bilinear estimate holds
\begin{eqnarray}\label{EstBlBasM1}
\nonumber||\partial_x(uv)||_{X^{-1/2+\delta,s_1-2\delta+\epsilon,s_2}}&\leq& CT^\mu\Big( ||u||_{X^{1/2,s_c^1,0}}
||v||_{X^{1/2,s_1,s_2}}\\
\nonumber&&\hspace{-2cm}+||u||_{X^{1/2,s_1,s_2}}||v||_{X^{1/2,s_c^1,0}}
+||u||_{X^{1/2,s_1,0}}||v||_{X^{1/2,s_c^1,s_2}}\\
&&+||u||_{X^{1/2,s_c^1,s_2}}||v||_{X^{1/2,s_1,0}}\Big)
\end{eqnarray}
this for some  $\delta>0$  small enough and $\epsilon>0$ such that $\epsilon<<\delta$.$\hfill{\Box}$
\end{proposition}
\section{Proof of Theorem \ref{Th1}}\label{RegEL}
\subsection{Exitence result}
Let $\phi\in H^{s_1,s_2}$ with  $s_1> -1/2$, $s_2\geq0$ and  $s_c^1\in ]-1/2,\min(0,s_1)]$. We suppose that $T\leq 1$, if  $u$ is a
solution of the integral equation $u=L(u)$ with
\begin{equation}\label{RegEL1}
L(u)=\psi(t)\Big[ W(t)\phi - \frac{\chi_{\mathbb R_+}(t)}{2} \int_0^t W(t-t')
\partial_x(\psi^2_T(t')u^2(t'))dt'\Big],
\end{equation}
then $u$ solve $KPB-II$- equation on $[0,T/2]$.We introduce the  Bourgain spaces defined by
 \begin{equation}\label{RegEL2}
Z_1=\{ u\in X^{1/2,s_1,s_2}; ||u||_{Z_1}=||u||_{ X^{1/2,s_1,0}}+
\gamma_1||u||_{ X^{1/2,s_1,s_2}}\},
\end{equation}
\begin{equation}\label{RegEL3}
Z_2=\{ u\in X^{1/2,s_1,0}; ||u||_{Z_2}=||u||_{ X^{1/2,s_c^1,0}}+
\gamma_2||u||_{ X^{1/2,s_1,0}}\},
\end{equation}
where
\begin{equation}\label{RegEL4}
\gamma_1=\frac{||\phi||_{H^{s_1,0}}}{||\phi||_{H^{s_1,s_2}}},\quad \gamma_2=\frac{||\phi||_
{H^{s_c^1,0}}}{||\phi||_{H^{s_1,0}}}.
\end{equation}
The goal to introduce two Bourgain spaces is to show in a first time that   there exists
$T_1=T(||\phi||_{H^{s_1,0}})$ and a solution $u$ of the equation (\ref{RegEL1}) in a ball of $Z_1$,
and then to solve   (\ref{RegEL1}) in $Z_2$ in order  to  check that the time of existence
$T=T(||\phi||_{H^{s_c^1,0}})$ with $s_c^1\in ]-1/2,0]$.\\
{\bf Step 1.} Resolution of (\ref{RegEL1}) in $Z_1$. By Proposition \ref{ProFree} and Proposition \ref{PrOp2}, it  results that,
\begin{equation}
||L(u)||_{X^{1/2,s_1,0}}\leq C||\phi||_{H^{s_1,0}}+
C||\partial_x(\psi_T^2(t)u^2)||_{X^{-1/2+\delta,s_1-2\delta+\epsilon,0}},
\end{equation}
\begin{equation}
||L(u)||_{X^{1/2,s_1,s_2}}\leq C||\phi||_{H^{s_1,s_2}}+
C||\partial_x(\psi_T^2(t)u^2)||_{X^{-1/2+\delta,s_1-2\delta+\epsilon,s_2}}.
\end{equation}
 By the proposition \ref{EstBlBas} and \ref{EstBlBasM}, we can deduce
\begin{equation}
||L(u)||_{X^{1/2,s_1,0}}\leq C||\phi||_{H^{s_1,0}}+
CT^{\mu}||\psi_T(t)u||^2_{X^{1/2,s_1,0}},
\end{equation}
\begin{equation}
||L(u)||_{X^{1/2,s_1,s_2}}\leq C||\phi||_{H^{s_1,s_2}}+
CT^{\mu}||\psi_T(t)u||_{X^{1/2,s_1,0}}||\psi_T(t)u||_{X^{1/2,s_1,s_2}},
\end{equation}
By Leibniz rule for fractional derivative and Sobolev inequalities in time  we have, for all $\epsilon >0$ and $0<T\leq1$, that
$$
||\psi_T(t)u||_{X^{1/2,s_1,s_2}}\leq C_\epsilon T^{-\epsilon}||u||_{X^{1/2,s_1,s_2}}.
$$
Taking $\epsilon=\mu/4$ we obtain,
 \begin{equation}\label{RegEL5}
||L(u)||_{X^{1/2,s_1,0}}\leq C||\phi||_{H^{s_1,0}}+
CT^{\mu/2}||u||^2_{X^{1/2,s_1,0}},
\end{equation}
\begin{equation}\label{RegEL6}
||L(u)||_{X^{1/2,s_1,s_2}}\leq C||\phi||_{H^{s_1,s_2}}+
CT^{\mu/2}||u||_{X^{1/2,s_1,0}}||u||_{X^{1/2,s_1,s_2}},
\end{equation}
By combining two estimates (\ref{RegEL5})  and (\ref{RegEL6})  we  obtain
\begin{equation}\label{RegEL7}
||L(u)||_{Z_1}\leq C(||\phi||_{H^{s_1,0}}+\gamma_1||\phi||_{H^{s_1,s_2}})+
CT^{\mu/2}||u||^2_{Z_1},
\end{equation}
Since $\partial_x(u^2)-\partial(v^2)=\partial_x[(u-v)(u+v)]$, in the same way we get
\begin{equation}\label{RegEL8}
||L(u)-L(v)||_{X_T^{1/2,s_1,0}}\leq CT^{\mu/2}||u-v||_{X^{1/2,s_1,0}}||u+v||_{X^{1/2,s_1,0}},
\end{equation}
\begin{eqnarray}\label{RegEL9}
\nonumber||L(u)-L(v)||_{X_T^{1/2,s_1,s_2}}&\leq& CT^{\mu/2}\Big(||u-v||_{X^{1/2,s_1,0}}
||u+v||_{X^{1/2,s_1,s_2}}\\
&+&
||u+v||_{X^{1/2,s_1,0}}||u-v||_{X^{1/2,s_1,s_2}}\Big),
\end{eqnarray}
it results that
\begin{equation}\label{RegEL10}
||L(u)-L(v)||_{Z_1}\leq CT^{\mu/2}||u-v||_{Z_1}||u+v||_{Z_1}.
\end{equation}
Hence, setting $T_1=(4C^2(||\phi||_{H^{s_1,0}}+\gamma_1||\phi||_{H^{s_1,s_2}}))^{-2/\mu}$
which yields by definition of $\gamma_1$ to  $T_1=(8C^2||\phi||_{H^{s_1,0}})^{-2/\mu}$ we deduce
from (\ref{RegEL7}) and  (\ref{RegEL10}) that $L$ is strictly contractive on the ball of radius
$r_1=2c(||\phi||_{H^{s_1,0}}+\gamma_1||\phi||_{H^{s_1,s_2}})$ in $Z_1$. This proves the existence of a unique solution $u_1$  to (\ref{RegEL1}) in  $X^{1/2,s_1,s_2}$ with $T_1=T(||\phi||_{H^{s_1,0}})$.\\
 Since $\phi \in H^{s_1,s_2}$, it follows that $\psi(.)W(.)\phi \in C([0,T_1], H^{s_1,s_2})$,
moreover  since $u_1\in X^{1/2,s_1,s_2}$, we can deduce from Proposition \ref{EstBlBas} that
$\partial_x(u_1^2)\in X_T^{-1/2+\delta,s_1-2\delta,s_2}$ and from (\ref{PrForReg1}) in Proposition
\ref{PrForReg} it results that
$$\int_0^tW(t-t')\partial_x(u_1^2(t'))dt'\in C([0,T_1], H^{s_1,s_2}).$$
Thus $u_1$ belongs $C([0,T_1], H^{s_1,s_2})$.\\
{\bf Step 2.} Resolution of (\ref{RegEL1}) in $Z_2$. Now proceeding exactly in the same way as above but in
$Z_2$, we obtain that $L$ is also strictly contractive on the  ball of radius
 $r_1=2c(||\phi||_{H^{s_c^1,0}}+\gamma_2||\phi||_{H^{s_1,0}})$ in $Z_2$  with
 $T_2=(4C^2(||\phi||_{H^{s_c^1,0}}+\gamma_2||\phi||_{H^{s_1,0}}))^{-1/\mu}$. Therefore by
definition of $\gamma_2$, it follows that $T_2=T(||\phi||_{H^{s_c^1,0}})$. Since obviously
 $H^{s_1,s_2}\subset H^{s_1,0}$, it follows that there exists a unique  solution $u_1$ to
(\ref{RegEL1}) in $C([0,T_2], H^{s_1,0})\cap X_T^{1/2,s_1,0}$ and $T_2=T(||\phi||_{H^{s_c^1,0}})$,
 $s_c^1\in ]-1/2,0]$. If we indicate by $T_*=T_{max}$ the maximum time of the
existence in $Z_1$ then   By uniqueness, we have $u_1=u_2$ on  $[0, \min(T_2,T_*)[$ and this gives
 that $T_*\geq T(||\phi||_{H^{s_c^1,0}})$.\\
 The continuity of map $\phi \longmapsto u$ from $H^{s_1,s_2}$ to  $X^{1/2,s_1,s_2}$ follows
from classical argument, while the continuity from $H^{s_1,s_2}$ to $C([0,T_1], H^{s_1,s_2})$ follows
again from Proposition \ref{PrForReg}.  The analyticity of the flow-map is a direct  consequence of the implicit function theorem. $\hfill{\Box}$
\subsection{Uniqueness.} The above contraction argument gives the uniqueness of the solution to the truncated integral equation (\ref{RegEL1}). We give here the argument of \cite{MoRi2} to deduce easily the uniqueness of the solution to the integral equation (\ref{FormDuh}).\\
 Let $u_1$, $u_2·\in X_T^{1/2,s_1,s_2}$ be two solution of the integral equation (\ref{FormDuh})
 on the time interval $[0,T]$ and let $\tilde u_1-\tilde u_2$ be an extension of  $u_1- u_2$ in
 $X^{1/2,s_1,s_2}$ such that
$$||\tilde u_1-\tilde u_2||_{X^{1/2,s_1,s_2}}\leq 2|| u_1- u_2||_{X_\gamma^{1/2,s_1,s_2}}$$
 with $0<\gamma\leq T/2$. It results by Proposition \ref{ProFree}
 and \ref{PrOp2} that,\\
$\hspace{-0,5 cm}||u_1-u_2||_{X_\gamma^{1/2,s_1,s_2}}$
\begin{eqnarray*}
&\leq& || \psi(t)\frac{\chi_{\mathbb R_+}(t)}{2} \int_0^t W(t-t')\partial_x\Big(\psi^2_\gamma(t')
\big(\tilde u_1(t')- \tilde u_2(t')\big)\big( u_1(t')+  u_2(t')\big) \Big)dt'
||_{X^{1/2,s_1,s_2}}\\
&\leq& C||\partial_x\Big(\psi_\gamma^2(t)
\big(\tilde u_1(t)- \tilde u_2(t)\big)\big( u_1(t)+  u_2(t)\big) \Big)||_{X^{-1/2+\delta,s_1-2\delta+\epsilon,s_2}}\\
&\leq& C\gamma^{\mu/2}||\tilde u_1- \tilde u_2||_{X^{1/2,s_1,s_2}} ||u_1+ u_2||_{X_T^{1/2,s_1,s_2}}
\end{eqnarray*}
for some $\mu >0$. Hence
$$||u_1-u_2||_{X^{1/2,s_1,s_2}_\gamma}\leq 2 C\gamma^{\mu/2} \Big(|| u_1||_{X^{1/2,s_1,s_2}_T}+
||u_2||_{X_T^{1/2,s_1,s_2}}\Big)||u_1-u_2||_{X^{1/2,s_1,s_2}_\gamma}.$$
Taking
 $\gamma \leq \Big( 4 C(|| u_1(t)||_{X^{1/2,s_1,s_2}_T}+
||u_2(t)||_{X_T^{1/2,s_1,s_2}})\Big)^{-\mu/2}$, this forces $u_1\equiv u_2$ on $[0,\gamma]$. Iterating this argument, one extends the uniqueness result on the whole time interval [0,T].$\hfill{\Box}$

\subsection{Global existence.}\
\\
By Proposition \ref{EstBlBasM}  $\partial_x(u^2)\in X^{-1/2+\delta,s_1+\epsilon-2\delta,s_2}$, therefore by Proposition \ref{PrForReg} we obtain that
$$\int_0^tW(t-t')\partial_x(u^2(t'))dt'\in C([0,T], H^{s_1+\epsilon,s_2}).$$
Note that $W(.)\phi\in  C(\mathbb R_+;H^{s_1,s_2})\cap C(\mathbb R_+^*;H^{\infty,s_2})$. Hence
$$u\in C([0,T];H^{s_1,s_2})\cap C(]0,T];H^{s_1+\epsilon,s_2}).$$
Recalling that $T=T(||\phi||_{H^{s_c^1,0}})$ with $s_1^c\in ]-1/2,\min(0,s_1)]$ and using the uniqueness result,
 we deduce by induction that $u\in C(]0,T];H^{\infty,s_2})$. This allows us to take the
 $L^2$-scalar product of (\ref{KPB}) with $u$, which shows that $t\longmapsto||u(t)||_{L^2}$ is
 nonincreasing on $]0,T]$. Since the time of local existence $T$ only depends on
$||\phi||_{H^{s_c^1,0}}$, this clearly  gives  that the solution is global in time.$\hfill{\Box}$
\section{ill- possedness results for the KPB-II equation}\label{ConExp}
In this section, we prove the ill- posedness result for the KPB-II equation  stated  in
 Theorem \ref{Th2}. We  start by   constructing a sequence of  initial data $(\phi_n)_n$ which
will   ensure  the non regularity of the map $\phi\longmapsto u$ from $H^{s,0}$ to $C([0,T],H^{s,0})$
for $s<-1/2$.
\subsection{Proof of Theorem \ref{Th2}}\label{ConExp1}
Let u be a solution of (\ref{KPB}) . Then we have
\begin{equation}\label{ConExp2}
u(\phi,t, x,y)=W(t)\phi(x,y) - \frac12 \int_0^t W(t-t')\partial_x(u^2(\phi,t',x,y))dt',
\end{equation}
suppose that the map is $C^2$. Since $u(0,t,x,y)=0$ it is easy to check  that
$$u_1(t,x,y)=:\frac{\partial u}{\partial \phi}(0,t,x,y)[h]=w(t)h$$
\begin{eqnarray}\label{ConExp3}
\nonumber  u_2(t,x,y)&:=&\frac{\partial^2 u}{\partial \phi^2}(0,t,x,y)[h,h]\\
&=&- \int_0^t W(t-t')u_1(t',x,y)\partial_x( u_1(t',x,y))dt'\\
\nonumber &=& - \int_0^t W(t-t')\partial_x(W(t')h)^2 dt'.
\end{eqnarray}
Due to the assumption of $C^2$-regularity of the map and since  that $u(0,t,x,y)=0$, we can write
 a formal Taylor expansion
\begin{equation}\label{ConExp4}
u(h,t,x,y)=u_1(t,x,y)[h] +u_2(t,x,y)[h,h]+ O(||h||_{H^{s,0}}^2),
\end{equation}
and we must  have
\begin{equation}\label{ConExp5}
||u_1(t,.,.)||_{H^{s,0}}\lesssim ||h||_{H^{s,0}},\quad
||u_2(t,.,.)||_{H^{s,0}}\lesssim ||h||^2_{H^{s,0}}.
\end{equation}
Taking the partial  Fourier transform with respect to $(x,y)$, it results that
\begin{eqnarray}\label{ConExp6}
\nonumber\mathcal F_{x\mapsto \xi, y\mapsto \eta}( u_2(t,.,.))&=&\int_0^t \exp[-|t-t'|\xi^2)
\exp(i(t-t')(\xi^3-\eta^2/\xi)]\\
&\times& (i\xi)\times\bigg [\mathcal F_{x\mapsto \xi, y\mapsto \eta}(u_1(t')\star u_1(t'))(\xi,\eta)\bigg]dt'
\end{eqnarray}
Note that\\
$$\hspace{-4cm} F_{x\mapsto \xi, y\mapsto \eta}\Big(u_1(t')\star u_1(t')\Big)(\xi,\eta)$$
\begin{eqnarray}\label{ConExp7}
\nonumber &=&F_{x\mapsto \xi, y\mapsto \eta}\Big(w(t')\phi\Big)\star
 F_{x\mapsto \xi, y\mapsto \eta}\Big(w(t')\phi\Big)\\
\nonumber&=& \int_{\mathbb R^2} \hat \phi(\xi-\xi_1,\eta-\eta_1)\hat\phi(\xi_1,\eta_1)
\exp\Big(-\big(\xi_1^2+(\xi-\xi_1)^2\big)t'\Big)\\
&&\hspace{1cm}\times \exp \Big(it\big(\xi_1^3 -\frac{\eta_1^2}{\xi_1} + (\xi-\xi_1)^2-
\frac{(\eta-\eta_1)^2}{\xi-\xi_1}\big)\Big)
d\xi_1 d\eta_1.
 \end{eqnarray}
  Now let $ P(\xi,\eta)=\xi^3-\eta^2/\xi$. Hence\\
$$\hspace{-8cm}\mathcal F_{x\mapsto \xi, y\mapsto \eta}( u_2(t,.,.))$$
\begin{eqnarray}\label{ConExp8}
\nonumber&=&\int_{\mathbb R^2} \hat \phi(\xi-\xi_1,\eta-\eta_1)\hat\phi(\xi_1,\eta_1)(i\xi)
e^{-t\xi^2}e^{itP(\xi,\eta)} \int_0^t e^{ it'\big(\xi_1^2 +(\xi-\xi_1)^2-\xi^2\big)}\\
&& \hspace{1cm}\times e^{it'\big(P(\xi_1,\eta_1)+P(\xi-\xi_1,\eta-\eta_1)-P(\xi,\eta)\big)}
dt'd\xi_1d\eta_1.
 \end{eqnarray}
Let $\chi(\xi,\xi_1,\eta,\eta_1)=P(\xi_1,\eta_1)+P(\xi-\xi_1,\eta-\eta_1)-P(\xi,\eta)$.
 A simple calulation show that
$$ \chi(\xi,\xi_1,\eta,\eta_1)=3\xi\xi_1(\xi-\xi_1)+\frac{(\eta\xi_1-\eta_1\xi)^2}
{\xi\xi_1(\xi-\xi_1)}.$$
therefore we can deduce that
\begin{eqnarray}\label{ConExp9}
\nonumber \mathcal F_{x\mapsto \xi, y\mapsto \eta}( u_2(t,.,.))&=&(i\xi)e^{itP(\xi,\eta)}
\int_{\mathbb R^2}\hat\phi(\xi_1,\eta_1)\hat \phi(\xi-\xi_1,\eta-\eta_1)\\
&&\times \frac{e^{-t(\xi_1^2+(\xi-\xi_1)^2)}e^{it\chi(\xi,\xi_1,\eta,\eta_1)}-e^{-\xi^2t}}
{-2\xi_1(\xi-\xi_1)+i\chi(\xi,\xi_1,\eta,\eta_1)}d\xi_1d\eta_1
\end{eqnarray}
it follows that
\begin{eqnarray}\label{ConExp10}
\nonumber ||u_2(t)||_{H^{s,0}}^2&:=& \int_{\mathbb R^2}(1+|\xi|^2)^s|\mathcal F_{x\mapsto \xi, y\mapsto \eta}( u_2(t,.,.)(\xi,\eta))|^2 d\xi d\eta\\
\nonumber&=&\int_{\mathbb R^2}|\xi|^2(1+|\xi|^2)^s\bigg|\int_{\mathbb R^2}\hat\phi(\xi_1,\eta_1)\hat \phi(\xi-\xi_1,\eta-\eta_1)\\
&&\times \frac{e^{-t(\xi_1^2+(\xi-\xi_1)^2)}e^{it\chi(\xi,\xi_1,\eta,\eta_1)}-e^{-\xi^2t}}
{-2\xi_1(\xi-\xi_1)+i\chi(\xi,\xi_1,\eta,\eta_1)}d\xi_1d\eta_1\bigg|^2 d\xi d\eta.
\end{eqnarray}
We chose now a real sequence of initial data $(\phi_N)_N$, $N>0$, defined through its Fourier
transform by
\begin{equation}\label{ConExp11}
\hat \phi_N(\xi,\eta)= N^{-3/2-s}(\chi_{D_{1,N}}+\chi_{D_{2,N}})
\end{equation}
where  $N$ is a positive  parameter  such that $N>>1$, and $D_{1,N}$, $D_{2,N}$ are the rectangles in
 $\mathbb R^2$ defined by
$$ D_{1,N}=[N/2,N]\times[-6N^2,6N^2],\quad D_{2,N}=[N,2N]\times [\sqrt{3}N^2,(\sqrt{3}+1)N^2].$$
It is simple to see that $||\phi_N||_{H^{s,0}}\sim 1$. Let us denote by $u_{2,N}$ the
 sequence  of the second iteration  $u_2$ associated with $\phi_N$. Setting
$$  K(t,\xi,\xi_1,\eta,\eta_1)= \frac{e^{-t(\xi_1^2+(\xi-\xi_1)^2)}
e^{it\chi(\xi,\xi_1,\eta,\eta_1)}-e^{-\xi^2t}}{-2\xi_1(\xi-\xi_1)+i\chi(\xi,\xi_1,\eta,\eta_1)},$$
 $||u_{2,N}(t)||_{H^{s,0}}^2$ can be split into three parts~:
$$||u_{2,N}(t)||_{H^{s,0}}^2= C(|f_1(t)|+|f_2(t)|+|f_3(t)|)$$
where
\begin{eqnarray*}
|f_1(t)|^{1/2}&=&CN^{-3-2s}\bigg[\int_{\mathbb R^2}|\xi|^2(1+|\xi|^2)^s\\
&&\times\bigg|\int_{\small{\begin{array}{l}
(\xi_1,\eta_1)\in D_{1,N}\\(\xi-\xi_1,\eta-\eta_1)\in D_{1,N}\end{array}}}
 K(t,\xi,\xi_1,\eta,\eta_1)d\xi_1 d\eta_1\bigg|^2 d\xi d\eta\bigg]^{1/2},
\end{eqnarray*}
\begin{eqnarray*}
|f_2(t)|^{1/2}&=&CN^{-3-2s}\bigg[\int_{\mathbb R^2}|\xi|^2(1+|\xi|^2)^s\\
&&\times\bigg|\int_{\small{\begin{array}{l}
(\xi_1,\eta_1)\in D_{2,N}\\(\xi-\xi_1,\eta-\eta_1)\in D_{2,N}\end{array}}}
 K(t,\xi,\xi_1,\eta,\eta_1)d\xi_1 d\eta_1\bigg|^2 d\xi d\eta\bigg]^{1/2},
\end{eqnarray*}
$$|f_3(t)|^{1/2}=CN^{-3-2s}\bigg[\int_{\mathbb R^2}|\xi|^2(1+|\xi|^2)^s\bigg|\int_{k(\xi,\eta)}
 K(t,\xi,\xi_1,\eta,\eta_1)d\xi_1 d\eta_1\bigg|^2 d\xi d\eta\bigg]^{1/2}.$$
where
\begin{eqnarray}\label{ConExp12}
\nonumber k(\xi,\eta)&=&\Big \{ (\xi_1,\eta_1): (\xi-\xi_1,\eta-\eta_1)\in D_{1,N},\space
 (\xi_1,\eta_1)\in D_{2,N}\Big\}\\&&
\nonumber\cup\Big \{ (\xi_1,\eta_1): (\xi_1,\eta_1)\in D_{1,N},\space
 (\xi-\xi_1,\eta-\eta_1)\in D_{2,N}\Big\}\\
&:=&k^1(\xi,\eta)\cup k^2(\xi,\eta).
\end{eqnarray}
Therefore, obviously
$$||u_{2,N}(t)||^2_{H^{s,0}}\geq C |f_3|.$$
Since we can write $\xi=\xi_1+(\xi-\xi_1)$, it follows that
\begin{eqnarray}\label{ConExp13}
\nonumber ||u_2(t)||_{H^{s,0}}^2&\geq&C N^{-4s-6}\int_{3N/2}^{3N}\int_{(\sqrt{3}-6)N^2}^{(\sqrt{3}+7)N^2}
|\xi|^2(1+|\xi|^2)^s\\
&&\times \bigg|\int_{k(\xi,\eta)} \frac{e^{-t(\xi_1^2+(\xi-\xi_1)^2)}e^{it\chi(\xi,\xi_1,\eta,\eta_1)}-e^{-\xi^2t}}
{-2\xi_1(\xi-\xi_1)+i\chi(\xi,\xi_1,\eta,\eta_1)}d\xi_1d\eta_1\bigg|^2 d\xi d\eta.
\end{eqnarray}
We  need to find a lower bound for the right-hand side  of (\ref{ConExp13}). Thus it is necessary to
 evaluate  the contribution of the function   $\chi(\xi,\xi_1,\eta,\eta_1)$  in $k(\xi,\eta)$.
This in the aim of the following lemma  which is inspired by \cite{MoSauTz}.
\begin{lemma}\label{AccFini}
Let $(\xi_1,\eta_1)\in k^1(\xi,\eta)$ or  $(\xi_1,\eta_1)\in k^2(\xi,\eta)$. For $N>>1$ we have
$$\big|\chi(\xi,\xi_1,\eta,\eta_1)\big|\lesssim N^3.$$
$\hfill{\Box}$
\end{lemma}
$\textbf{Proof of lemma \ref{AccFini}.}$
By  definition of the fonction $\chi(\xi,\xi_1,\eta,\eta_1)$ we can write
\begin{equation}\label{AccFini1}
\big|\chi(\xi,\xi_1,\eta,\eta_1)\big|\leq \big|\chi_1(\xi,\xi_1,\eta,\eta_1)\big|+
\big|6\xi\xi_1(\xi-\xi_1)\big|,
\end{equation}
where $$ \chi_1(\xi,\xi_1,\eta,\eta_1)=3\xi\xi_1(\xi-\xi_1) -\frac{(\eta\xi_1-\eta_1\xi)^2}
{\xi\xi_1(\xi-\xi_1)}.$$
Now let $(\xi_1,\eta_1)\in k^1(\xi,\eta)$  i.e. $(\xi-\xi_1,\eta-\eta_1)\in D_{1,N}$ and $(\xi_1,\eta_1)\in D_{2,N}$.\\
Let $\xi\in \mathbb R$ such that $(\xi-\xi_1)\in [N/2,N]$ and we fix  $(\xi_1,\eta_1)\in D_{2,N}$.
we will seek a $\eta^*(\xi,\xi_1,\eta_1)$ such that
$\chi_1(\xi,\xi_1,\eta^*(\xi,\xi_1,\eta_1),\eta_1)=0$ and
$\big|\eta^*(\xi,\xi_1,\eta_1)-\eta_1\big|\leq 6N^2$. Indeed, we choose
$$\eta^*(\xi,\xi_1,\eta_1)=\eta_1 +\frac{(\xi-\xi_1)(\eta_1-\sqrt{3}\xi\xi_1)}{\xi_1}.$$
And thus
$$
\big|\eta^*(\xi,\xi_1,\eta_1)-\eta_1\big|\leq \frac{|\xi-\xi_1|}{|\xi_1|}
\big|\eta_1-\sqrt{3}\xi_1^2-\sqrt{3}\xi_1(\xi-\xi_1)\big|
$$
We recall that $\eta_1  \in [\sqrt{3}N^2,(\sqrt{3}+1)N^2]$ and $\xi_1\in[N,2N]$. Therefore
 It follows that
$$\sqrt{3} \xi_1^2\in [\sqrt{3}N^2,4\sqrt{3}N^2]$$
and, we have
$$
\big|\eta_1-\sqrt{3}N^2\big|\leq 3\sqrt{3} N^2.
$$
 Since $|\xi_1|\leq 2N$ and $|\xi-\xi_1|\geq N/2$, it results that
$$\big|\eta^*(\xi,\xi_1,\eta_1)-\eta_1\big|\leq1/4\Big(3\sqrt{3}N^2 +2\sqrt{3}N^2\Big)\leq 6N^2.$$
Now by the mean value theorem we can write
$$\chi_1(\xi,\xi_1,\eta,\eta_1)=\chi_1(\xi,\xi_1,\eta^*(\xi,\xi_1,\eta_1),\eta_1)+
\big(\eta-\eta^*(\xi,\xi_1,\eta_1)\big)\frac{\partial \chi_1}{\partial\eta}
(\xi,\xi_1,\bar{\eta},\eta_1)$$
where $\bar{\eta}\in[\eta,\eta^*(\xi,\xi_1,\eta_1)]$. Therefore
$$\big|\chi_1(\xi,\xi_1,\eta,\eta_1)\big|=
\big|\eta-\eta^*(\xi,\xi_1,\eta_1)\big|\Big|\frac{2\xi_1(\bar{\eta}\xi_1-\eta_1\xi)}
{\xi\xi_1(\xi-\xi_1)}\Big|.$$

Since $\big|\eta-\eta^*(\xi,\xi_1,\eta_1)\big|\leq \big|\eta-\eta_1\big| +
\big|\eta_1-\eta^*(\xi,\xi_1,\eta_1)\big|\leq CN^2$, it follows that
\begin{eqnarray*}
\big|\chi_1(\xi,\xi_1,\eta,\eta_1)\big|&\lesssim&
|\xi_1|\big|\eta-\eta^*(\xi,\xi_1,\eta_1)\big|\Big|\frac{(\bar{\eta}-\eta_1)\xi_1-\eta_1(\xi-\xi_1)}
{\xi\xi_1(\xi-\xi_1)}\Big|\\
&\lesssim&
N^3\bigg(\frac{|(\bar{\eta}-\eta_1)\xi_1|}
{|\xi\xi_1(\xi-\xi_1)|}+\frac{|\eta_1(\xi-\xi_1)|}
{|\xi\xi_1(\xi-\xi_1)|}\bigg)\\
&\lesssim& N^3\bigg(\frac{(\sqrt{3}+1)N^2}{N^2} +C\frac{N^2}{N^2}\bigg)\\
&\lesssim& N^3
\end{eqnarray*}
by the relation of (\ref{AccFini1}) it results that
$$\big|\chi(\xi,\xi_1,\eta,\eta_1)\big|\leq CN^3+ 6|\xi||\xi_1||\xi-\xi_1|,$$
since one has $|\xi|\sim|\xi_1|\sim|\xi-\xi_1|\sim N$, it follows that
$$\big|\chi(\xi,\xi_1,\eta,\eta_1)\big|\leq CN^3.$$
Now, in the other case where $(\xi_1,\eta_1)\in k^2(\xi,\eta)$  i.e.  $(\xi_1,\eta_1)\in D_{1,N}$ and $(\xi-\xi_1,\eta-\eta_1)\in D_{2,N}$,
follows from first case since we can write
 $(\xi_1,\eta_1)=(\xi-(\xi-\xi_1),\eta-(\eta-\eta_1))\in D_{1,N}$ and  that
$$\chi_1(\xi,\xi_1,\eta,\eta_1)=\chi_1(\xi,\xi-\xi_1,\eta,\eta-\eta_1).$$
This completes the proof of the Lemma.$\hfill{\Box}$\\
Let us now end the proof of  Theorem \ref{Th2}. Note that for any $\xi\in [3N/2,3N]$ and
$[(\sqrt{3}-6)N^2,(\sqrt{3}+7)N^2]$, we have mes$\big(K(\xi,\eta)\big)\geq CN^3$. We recall that we
have
\begin{eqnarray}\label{ConExp14}
\nonumber ||u_{2,N}(t)||_{H^{s,0}}^2&\geq&C N^{-4s-6}\int_{3N/2}^{3N}
\int_{(\sqrt{3}-6)N^2}^{(\sqrt{3}+7)N^2}|\xi|^2(1+|\xi|^2)^s\\
&\times& \bigg|\int_{k(\xi,\eta)}
 \frac{e^{-\xi^2t}\Big[e^{\big(-2\xi_1(\xi-\xi_1)+it\chi(\xi,\xi_1,\eta,\eta_1)\big)}-1\Big]}
{-2\xi_1(\xi-\xi_1)+i\chi(\xi,\xi_1,\eta,\eta_1)}d\xi_1d\eta_1\bigg|^2 d\xi d\eta.
\end{eqnarray}
Now, we choose a  sequence of  times  $(t_N)_N$ defined  by
$$t_N=\frac{1}{N^{3+\epsilon_0}}, \quad 0<\epsilon_0<<1 \text{ }(\text{fixed}).$$
 For $N>>1$ it is clear
\begin{equation}\label{ConExp15}
e^{-\xi^2t_N}\sim e^{-N^2t_N}\sim e^{-\frac{1}{N^{1+\epsilon_0}}}>C.
\end{equation}
Moreover, by  Lemma \ref{AccFini}, it follows that
 $\big|-2\xi_1(\xi-\xi_1)+it\chi(\xi,\xi_1,\eta,\eta_1)\big|\leq N^2+N^3\leq CN^3$. Hence
\begin{eqnarray}\label{ConExp16}
\nonumber\bigg|\frac{e^{\big(-2\xi_1(\xi-\xi_1)+it\chi(\xi,\xi_1,\eta,\eta_1)\big)}-1}
{-2\xi_1(\xi-\xi_1)+i\chi(\xi,\xi_1,\eta,\eta_1)}\bigg|&=& |t_N| +O(|t_N|^2N^3)\\
&=&\frac{1}{N^{3+\epsilon_0}} +O(\frac{1}{N^{3+2\epsilon_0}})
\end{eqnarray}
By combining the relations (\ref{ConExp15}), (\ref{ConExp16}), we obtain
$$
\nonumber\bigg|\int_{k(\xi,\eta)} \frac{e^{-\xi^2t}\Big[e^{\big(-2\xi_1(\xi-\xi_1)+it
\chi(\xi,\xi_1,\eta,\eta_1)\big)}-1\Big]}
{-2\xi_1(\xi-\xi_1)+i\chi(\xi,\xi_1,\eta,\eta_1)}d\xi_1d\eta_1\bigg|
$$
\begin{equation}\label{ConExp17}
\geq C\text{mes}\big(k(\xi,\eta)\big)\times\frac{1}{N^{3+\epsilon_0}}\geq CN^{-\epsilon_0}.
\end{equation}
By virtue of (\ref{ConExp14}), it results that
\begin{eqnarray}
\nonumber ||u_{2,N}(t_N)||_{H^{s,0}}^2&\geq&C N^{-4s-6}\int_{3N/2}^{3N}
\int_{(\sqrt{3}-6)N^2}^{(\sqrt{3}+7)N^2}|\xi|^2(1+|\xi|^2)^s d\xi d\eta \times N^{-2\epsilon_0}\\
\nonumber&\geq& C N^{-6-4s}N^{2s}N^2N^3N^{-2\epsilon_0}\\
\nonumber&\geq&  CN^{-1-2\epsilon_0-2s}
\end{eqnarray}
and, hence
$$1\sim ||\phi_N||_{H^{s,0}}^2\geq ||u_{2,N}(t_N)||^2_{H^{s,0}}\geq N^{-1-2\epsilon_0-2s}.$$
This  leads to  a contradiction for $N>>1$, since we have   $-1-2\epsilon_0-2s>0$
for $s\leq -1/2+\epsilon_0$. This completes the proof of Theorem \ref{Th2}. $\hfill{\Box}$\\\\
$\textbf{Acknowledgments.}$ I  would like  to thank Luc Molinet for his encouragement,
advice, help and for the rigorous attention to this paper.

\vspace{1cm}
\author{Bassam Kojok~:}
\end{document}